\documentclass[11pt]{amsart}
\usepackage{graphicx}
\usepackage{latexsym,amsmath,amsfonts,amscd, amsthm}
\usepackage{color}
\usepackage[colorlinks=true]{hyperref}
\hypersetup{urlcolor=blue, citecolor=red}
\usepackage{tikz}
\usetikzlibrary{arrows,backgrounds,snakes}
\usepackage{epsfig}
\usepackage{dsfont}
\usepackage{amssymb}
\usepackage{xcolor}
\usepackage{comment}
\usepackage{subfigure}
\usepackage{multirow}

\newtheoremstyle{plainNoItalics}{}{}{\normalfont}{}{\bfseries}{.}{ }{}
\theoremstyle{plain}
\newtheorem{thm}{Theorem}[section]

\newtheorem{defn}[thm]{Definition}
\newtheorem{rem}[thm]{Remark}
\newtheorem{prop}[thm]{Proposition}

\newcommand{\f}{\frac}
\newcommand{\ra}{\rightarrow}
\newcommand{\beq}{\begin{equation}}
\newcommand{\eeq}{\end{equation}}
\newcommand{\beqa}{\begin{eqnarray}}
\newcommand{\eeqa}{\end{eqnarray}}
\newcommand{\bit}{\begin{itemize}}
\newcommand{\eit}{\end{itemize}}
\newcommand{\bedef}{\begin{defn}}
\newcommand{\edefn}{\end{defn}}
\newcommand{\bpro}{\begin{prop}}
\newcommand{\epro}{\end{prop}}

\newcommand{\divx}{{\nabla\cdot}}

\newcommand{\df}{\partial}

\newcommand{\mD}{{\mathcal D}}

\newcommand{\mO}{{\mathcal O}}

\newcommand{\mI}{{\mathbb  I}}

\newcommand\vave[1]{{\langle{#1}\rangle}}

\newcommand{\eps}{\varepsilon}

\newcommand{\vr}{{\vec{r}}}
\newcommand{\vO}{{\vec{\Omega}}}

\newcommand{\xL}{{x_{i-\frac{1}{2}}}}
\newcommand{\xR}{{x_{i+\frac{1}{2}}}}
\newcommand{\yL}{{y_{j-\frac{1}{2}}}}
\newcommand{\yR}{{y_{j+\frac{1}{2}}}}

\newcommand{\testR}{{\alpha}}   
\newcommand{\testP}{{\bf Z}}   
\newcommand{\testG}{{\beta}}
\newcommand{\testE}{{\gamma}}   

\newcommand{\pI}{\partial I}
\newcommand{\pT}{\partial T}
\newcommand{\pU}{\partial U}
\newcommand{\pt}{\partial t}

\newcommand{\tC}{\tilde{C}_v}

\newcommand{\bk}{{\bf k}}

\newcommand{\bn}{{\bf n}}
\newcommand{\bq}{{\bf q}}

\newcommand{\ssig}{{\sqrt{\sigma_0}}}
\newcommand{\rsig}{{\tilde{\sigma}}}
\newcommand{\reps}{{\tilde{\varepsilon}}}

\email{txiong@xmu.edu.cn (T. Xiong)}
\email{sun\underline{ }wenjun@iapcm.ac.cn (W. Sun)}
\email{shi\underline{ }yi@iapcm.ac.cn (Y. Shi)}
\email{song\underline{ }peng@iapcm.ac.cn (P. Song)}

\title[High order AP DG for gray RTE]{High order asymptotic preserving discontinuous Galerkin methods for gray radiative transfer equations} 


\begin{document}
	
\maketitle

\medskip
\centerline{\scshape Tao Xiong\footnote{the corresponding author.}}
\medskip
{\footnotesize
\centerline{School of Mathematical Sciences, Fujian Provincial Key Laboratory of Mathematical Modeling} 
\centerline{and High-Performance Scientific Computing, Xiamen University}
\centerline{Xiamen, Fujian 361005, P.R. China}
}

\bigskip

\centerline{\scshape Wenjun Sun \quad Yi Shi \quad Peng Song}
\medskip
{\footnotesize
	\centerline{Institute of Applied Physics and Computational Mathematics} 
	\centerline{P.O. Box 8009, Beijing 100088, P.R. China}
}

\bigskip

\begin{abstract}
In this paper, we will develop a class of high order asymptotic preserving (AP) discontinuous Galerkin (DG) methods for nonlinear time-dependent gray radiative transfer equations (GRTEs). Inspired by the work \cite{Peng2020stability}, in which stability enhanced high order AP DG methods are proposed for linear transport equations, we propose to pernalize the nonlinear GRTEs under the micro-macro decomposition framework by adding a weighted linear diffusive term. In the diffusive limit, a hyperbolic, namely 
$\Delta t=\mO(h)$ where $\Delta t$ and $h$ are the time step and mesh size respectively, instead of parabolic  
$\Delta t=\mO(h^2)$ time step restriction is obtained, which is also free from the photon mean free path. The main new ingredient is that we further employ a Picard iteration with a predictor-corrector procedure, to decouple the resulting global nonlinear system to a linear system with local nonlinear algebraic equations from an outer iterative loop. Our scheme is shown to be asymptotic preserving and asymptotically accurate. Numerical tests for one and two spatial dimensional problems are performed to demonstrate that our scheme is of high order, effective and efficient.
\end{abstract}

\vspace{0.1cm}

\tableofcontents

\section{Introduction}
\label{sec1}
\setcounter{equation}{0}
\setcounter{figure}{0}
\setcounter{table}{0}

Radiative transfer equations (RTEs) are a type of kinetic scale modeling equations, which are used to describe the time evolution of radiative intensity and energy transfer of a radiation field with its background material \cite{chandrasekhar,pomraning1973radiation}. The system has many applications in astrophysics, inertial confinement fusion (ICF), plasma physics and so on. It has attracted a lot of attention for numerical studies due to its importance but high complexity. 

Gray radiative transfer equations (GRTEs) are a type of simplified RTEs for gray photons and coupled to the background with the material temperature. Due to its high dimensionality and the photons are traveling in the speed of light, a popular numerical method for simulating the GRTEs in literature is the implicit Monte Carlo method, see \cite{Fleck1971,Gentile2001,McClarren2009,Densmore2011,Shi2018,Shi2020} and references there in. On the other hand, to avoid numerical noises from the Monte Carlo method, deterministic methods are also developed. However, several difficulties arise. Firstly, for the GRTEs with a high or thick opacity background material, there are severe interactions between radiation and material, in which case, the photon mean free path is approaching zero and the diffusive radiation behavior becomes to dominate. Numerical methods for streaming transport equations in the low or thin opacity material will suffer from very small time step restriction in this diffusive regime. Secondly, the photons are traveling in a speed of light, so that the system also deserves implicit treatments for the transport term. Thirdly, the time evolution of the radiative intensity is coupled with the background material temperature. Its local thermal equilibrium is a Planckian, which is a nonlinear function of the material temperature. In the diffusive limit as the photon mean free path vanishing, it converges to a nonlinear diffusive equation. Due to the high dimensionality and strong coupled nonlinearity, fully implicit schemes are very difficult to solve, which require huge computational cost unaffordable even on modern computers.  

An efficient approach to deal with these difficulties is the asymptotic preserving (AP) scheme. AP schemes were first studied in the numerical solution of steady neutron transport problems \cite{larsen1987asymptotic,larsen1989asymptotic,jin1991,jin1993}, and later applied to nonstationary transport problems \cite{klar1998asymptotic,jin2000uniformly}. For a scheme to be AP, in the diffusive scaling as we considered here, it means that when the spatial mesh size $h$ and time step $\Delta t$ are fixed, as the photon mean free path approaching zero, the scheme automatically becomes a consistent discretization for the limiting diffusive equation, see \cite{jin1999efficient,jin2010asymptotic} for more concrete definitions of AP schemes for multi-scale kinetic equations. Recently AP schemes are further developed under different scenorios. For linear transport equations in the diffusive scaling, a class of AP schemes coupled with discontinuous Galerkin discretizations in space and globally stiffly accurate implicit-explicit (IMEX) schemes were developed in \cite{Jang2014analysis,Jang2015high,Peng2020stability,Peng2020asymptotic} under the micro-macro decomposition framework \cite{lemou2010new}, while in space using finite difference discretizations under the same framework was developed in \cite{Cory2019}. Crestetto et. al. recently proposed to decompose the solution of the linear transport equation in the diffusive regime in the time evoluation of an equilibrium state plus an perturbation, and a Monte Carlo solver for the perturbation with an Eulerian method for the equilibrium part was designed \cite{Crestetto2019}, which is asymptotically complexity diminishing. AP schemes are also coupled with gas-kinetic schemes for linear kinetic model \cite{Mieussens2013} and nonlinear RTEs \cite{Sun2015an,Sun2015frequency,Sun2017multidimensional,Shi2020an,Xu2020}, etc. 
However, to the best of our knowledge, high order schemes for time-dependent nonlinear RTEs, which is AP for the vanishing photon mean free path, is rare in literature, except \cite{Hammer2019} which is based on multiscale high order/low order (HOLO) algorithms \cite{Chacon2019} and \cite{Maginot2016} which is based on high order S-stable diagonally implicit Runge-Kutta (SDIRK) methods and linearization on the Planck function, as some of such few examples.    

In this paper, inspired by the work \cite{Peng2020stability,Jang2015high}, we try to develop a class of high order AP DG-IMEX schemes for the GRTEs.
Here in space, we use DG discretizations. The first DG scheme is developed for neutron transport equation by Reed and Hill \cite{Reed1973}. Later DG has been widely used to study stationary and non-stationary transport equations, e.g. \cite{Adams2001,guermond2010asymptotic,Jang2015high,Yuan2016,Zhang2019,Cory2019,Peng2020stability,Cory2020,Guermond2020,Zhang2020high}. In time,  globally stiffly accurate IMEX schemes \cite{ascher1997implicit} are used.
The same space and time discretizations have been proposed for linear transport equations in the diffusive scaling in \cite{Jang2015high} and analized in \cite{Jang2014analysis}. Recently, Peng et. al. \cite{Peng2020stability} proposed a class of stability enhanced schemes by adding a weighted diffusive term on both sides but numerically discretizing them differently in time. The schemes have been analized in \cite{Peng2020stability2}, which show that the stability enhanced approach has a hyperbolic time step condition, namely $\Delta t = \mO(h)$ where $\Delta t$ and $h$ are the time step size and mesh size respectively, instead of a parabolic time step condition $\Delta t=\mO(h^2)$ in \cite{Jang2014analysis}. Further, to avoid the ad-hoc weighted diffusive term, Peng and Li proposed a class of AP IMEX-DG-S schemes based on the Schur complement \cite{Peng2020asymptotic}. Here, the stability-enhanced approach is borrowed, which is more straightforward than the scheme based on the Schur complement, when applied to the nonlinear GRTEs. The new challenge for GRTEs here is that the radiative intensity is nonlinearly coupled to the background material temperature. To effectively and efficiently deal with the nonlinearity appearred in the Planck function, we use a Picard iteration with a predictor-corrector procedure as is used in \cite{tang2020} for the capturing the right front propogation for GRTEs. The global nonlinear system is decoupled to a linear system with nonlinear algebraic equations only in each element from an outer iterative loop. We will formally prove that our scheme for GRTEs is AP and also asymptotically accurate (AA). Numerically we will show it is high order in space and in time, and can effectively and efficiently solve the nonlinear GRTEs.

The rest of the paper is organized as follows. In Section 2, the GRTEs are revisited and reformulated by a micro-macro decomposition.
DG space and IMEX time discretizations are described in Section 3. Formal AP and AA analyses are given in Section 4. We will perform some numerical tests in Section 5, followed by a conclusion in Section 6.

%
\section{Model equation}
\label{sec2}
\setcounter{equation}{0}
\setcounter{figure}{0}
\setcounter{table}{0}

\subsection{Gray radiative transfer equation}
\label{sec2.1}

Here we consider the gray radiative transfer equations in the scaled form \cite{Sun2015an, pomraning1973radiation}
\beq
\label{rte}
\left\{
\begin{array}{l}
\f{\eps^2}{c}\f{\pI}{\pt}+\eps\,\vec{\Omega}\cdot\nabla I =\sigma\left(\Phi-I\right), \\ \, \\
\eps^2\,C_v \f{\pT}{\pt}\equiv \eps^2\f{\pU}{\pt}=\sigma \, |\vO| \,\left(\vave{I}-\Phi\right),
\end{array}
\right.
\eeq
This system describes the radiative transfer and energy exchange between the radiation and the material. Here $t$ is the time variable, $\vr$ is the position variable and $\vO$ is the direction of traveling of the photons. $I(t,\vr,\vO)$ is the radiative intensity in the direction of $\vO$, $T(t,\vr)$ is the material temperature, $\sigma(t,\vr)$ is the opacity, $C_v(t,\vr)$ is the scaled heat capacity and $U(t,\vr)$ is the material energy density. $\eps > 0$ is the Knudsen number defined as the ratio of the photon mean free path over the characteristic length of space. We denote $\langle \cdot \rangle$ as the integral average over the angular $\vO$, namely
\beq
\vave{\cdot} = \frac{1}{|\vO|}\int \cdot \,d\vO,
\eeq
and $\Phi=acT^4/|\vO|$, where $a$ is the radiation constant, $c$ is the scaled speed of light. For simplicity, the internal source and scattering are omitted. 

The spatial position vector $\vr$ is usually presented by the Cartesian coordinate with $\vr=(x,y,z)$, while the direction $\vO=(\zeta,\eta,\mu)$ is described by a polar angle $\theta$ measured with respect to any fixed direction in space (e.g., the $z$ axis) and a corresponding azimuthal angle $\phi$, then
\[
\mu = \cos\theta, \quad \zeta=\sin\theta \cos\phi, \quad \eta = \sin\theta \sin\phi,
\]
and
\[
d\vr=dx\,dy\,dz,  \quad
d\vO=\sin\theta\,d\theta\,d\phi=-d\mu\,d\phi,
\]
In the one-dimensional case, \eqref{rte} reduces to
\beq
\label{rte1d}
\left\{
\begin{array}{l}
	\f{\eps^2}{c}\f{\pI}{\pt}+\eps\,\mu \, I_x =\sigma\left(\Phi-I\right), \\ \, \\
	\eps^2 C_v \f{\pT}{\pt}\equiv \eps^2\f{\pU}{\pt}=2\sigma\left(\vave{I}-\Phi\right),
\end{array}
\right.
\eeq
with $\mu\in[-1,1]$. While in the two dimensional case, it becomes
\beq
\label{rte2d}
\left\{
\begin{array}{l}
	\f{\eps^2}{c}\f{\pI}{\pt}+\eps\,(\zeta \, I_x + \eta \, I_y)=\sigma\left(\Phi-I\right), \\ \, \\
	\eps^2 C_v \f{\pT}{\pt}\equiv \eps^2\f{\pU}{\pt}=2\pi\,\sigma\left(\vave{I}-\Phi\right),
\end{array}
\right.
\eeq
with
\[
\zeta=\sqrt{1-\mu^2}\,\cos\phi\in[-1,1], \quad \eta=\sqrt{1-\mu^2}\,\sin\phi\in[-1,1],\quad \mu \in [0,1],\quad \phi\in[0,2\pi].
\]

\eqref{rte} is a relaxation model for the radiative intensity to its local thermodynamic equilibrium, at where the emission source is a Planckian at the local material temperature $\sigma\f{acT^4}{|\vO|}$. The material temperature $T(t,\vr)$ is related to the material energy density $U(t,\vr)$ by
\beq
\f{\pU}{\pT}=C_v > 0,
\eeq
where $C_v$ is the heat capacity.

It has been shown in \cite{Larsen1983} that as the Knudsen number $\eps \ra 0$, away from boundaries and initial times, the intensity $I$ will approach to the Planckian at the local temperature equilibrium, 
\beq
I^{(0)}=\f{1}{|\vO|}ac(T^{(0)})^4,
\eeq
and the corresponding local temperature $T^{(0)}$ satisfies the following nonlinear diffusive equation
\beq
\label{diffeq}
\f{\partial }{\pt}U(T^{(0)})+
\f{\partial}{\pt}a(T^{(0)})^4=\nabla\cdot\f{1}{3\sigma}\nabla ac(T^{(0)})^4.
\eeq

\subsection{Micro-macro decomposition}
We consider to reformulate the radiative transfer equations \eqref{rte} by a micro-macro decomposition \cite{lemou2010new, Jang2015high}. We start to split the radiative intensity $I$ as
\beq
\label{mmdI}
I = \rho + \frac{\eps}{\ssig} \, g,
\eeq
where $\rho = \Pi\, I=\langle I \rangle$, $\eps\,g/\ssig=(\mI-\Pi)I$ is the perturbation. $\mI$ is the identity and $\Pi$ is an orthogonal projection. $\sigma_0$ is a constant which is defined as a referred opacity.
$\ssig$ is added in the consideration that when $\sigma$ is relatively large, \eqref{rte} will also be close to the thermodynamic equilibrium, which is different from \cite{Jang2015high,Peng2020stability}. If we integrate over the travel direction $\vO$ of photons, subtracting it from \eqref{rte}, we will obtain the following micro-macro decomposition system
\beq
\label{eq:mmd}
\left\{
\begin{array}{l}
	\f{\partial \rho}{\pt}+\frac{c}{\ssig}\,\nabla\cdot \langle \vO\, g \rangle =\f{c\,\rsig}{\eps^2/\sigma_0}\left(\Phi-\rho\right), \\ \, \\
	\frac{\eps^2}{\sigma_0} C_v \f{\pT}{\pt}=\rsig \, |\vO|(\rho-\Phi), \\ \, \\
	g_t+\f{c}{\eps}(\mI-\Pi)\nabla\cdot(\vO\, g)+\f{c\,\ssig}{\eps^2}\nabla\cdot(\vO\,\rho)=-\f{c \,\rsig}{\eps^2/\sigma_0}g,
\end{array}
\right.
\eeq
where $\rsig=\sigma/\sigma_0$ is a relative opacity compared to the reference opacity $\sigma_0$ , and $\sigma_0$ will be specified later. 

In \eqref{eq:mmd}, if we assume all variables, such as $\rho$ and $g$, are of $\mathcal{O}(1)$, formally if $\eps/\ssig\ra 0$, in the leading order we get $\rho=\Phi$ and $g=-\f{\ssig}{\sigma}\nabla\cdot(\vO\,\rho)$ from the second and third equations respectively. Substituting them into
the first equation, replacing the right hand side by $-c\,\frac{C_v}{|\vO|}\f{\pT}{\pt}$ from the second equation, the nonlinear diffusive equation \eqref{diffeq} is directly implied, which shows the advantage of this micro-macro reformulation. A more rigorous asymptotic analysis will be given in the Section \ref{sec3.5}. For more discussions on micro-macro decomposition, see \cite{lemou2010new,Jang2015high}. 

In \cite{Jang2014analysis, Jang2015high}, a family of high order asymptotic preserving (AP) schemes was proposed and analyzed for the linear transport equations (corresponding to $\Phi=\rho$ in the second equation of \eqref{eq:mmd}, namely, the temperature is at the equilibrium which is time independent). This type of scheme is based on nodal discontinuous Galerkin (DG) spatial discretization and globally stiffly accurate implicit-explicit (IMEX) temporal discretization. However, the scheme as $\eps\ra 0$ in the diffusive limit is intrinsically an explicit time discretization for the limiting heat equation and numerical stability requires a time step size $\Delta t\sim h^2$ when $\eps\ll1$, where $h$ is the spatial mesh size. This time step condition is quite stringent for the consideration of computational efficiency. Recently, Peng et. al. \cite{Peng2020stability,Peng2020stability2} proposed a stability enhanced approach by adding a weighted diffusive term, based on the knowledge of the diffusive limit, which can improve the time step condition up to $\Delta t\sim h$ for stability and the scheme becomes much more efficient.

Now let us follow the idea in \cite{Peng2020stability} to further reformulate \eqref{eq:mmd} by adding a linear diffusive term on both sides in the first equation, with the referred opacity $\sigma_0$
\beq
\label{eq:mmd2}
\left\{
\begin{array}{l}
	\f{\partial \rho}{\pt}+c\,\left(\frac{1}{\ssig}\nabla\cdot \langle \vO\, g \rangle +\f{\omega}{3\sigma_0}\Delta \rho\right)=\f{c\,\rsig}{\reps^2}\left(\Phi-\rho\right)+c\,\f{\omega}{3
		\sigma_0}\Delta \rho, \\ \, \\
	\reps^2 \tC \f{\pT}{\pt}=\rsig (\rho-\Phi), \\ \, \\
	g_t+\f{c/\ssig}{\reps}(\mI-\Pi)\nabla\cdot(\vO\, g)+\f{c/\ssig}{\reps^2}\nabla\cdot(\vO\,\rho)=-\f{c\,\rsig}{\reps^2}g,
\end{array}
\right.
\eeq
where $\tC=C_v/|\vO|$ and $\reps=\eps/\ssig$ are taken for shorthand notation. $\omega$ is a non-negative weighted function of $\reps$ and it is bounded and independent of $\vr$. It suggests to take $\omega=\exp(-\eps/h)$ or $\omega=\exp(-\eps^2/\Delta t)$ in \cite{Peng2020stability}. In this work, we take $\omega=\exp(-\reps/h)$. The referred value $\sigma_0$ can be seen as a linear diffusive pernalization coefficient, e.g., as suggested in \cite{wang2020}, we may take $1/\sigma_0\ge 0.54 \max(1/\sigma)$. The two added linear diffusive terms will be numerically treated differently, one is explicit and the other is implicit.

Before we continue to numerical discretizations for \eqref{eq:mmd2}, we first mention the angular approximation. For the angular direction $\vO$, we simply use the discrete ordinate method, also known as $S_N$ method \cite{Lathrop1968}. In the one dimensional case, $\vO$ simply becomes $\mu\in[-1,1]$, a Gauss-Legendre quadrature rule with weights $w_m$ and nodes $\mu_m$ is used, for $m=1,\ldots,N_{g}$. In the two dimensional case, a Legendre-Chebyshev quadrature rule with weights $w_m$ and nodes $\vO_m=(\zeta_m,\eta_m)$, for $m=1,\ldots,N_g=N_lN_c$ is used. The nodes $\vO_m=(\zeta_m,\eta_m)$ are given by
\[
\zeta_m=\sqrt{1-\mu^2_i}\,\cos\phi_j, \quad \eta_m=\sqrt{1-\mu^2_i}\,\sin\phi_j, \quad m=(i-1)N_c+j,
\]
where $\mu_i$, $i=1,\ldots,N_l$, denote the roots of the Legendre polynomial of degree $N_l$. 
$\phi_j=(2j-1)\pi/N_c$, $j=1,\ldots,N_c$ are the nodes based on a Chebyshev polynomial.
In the following, for easy presentation, we still keep $\vO$ continuous and focus on the numerical discretizations in time and space.

\begin{rem}
	We notice that taking the micro-macro decomposition in \eqref{mmdI}, which is different from \cite{Jang2015high,Peng2020stability}, it has the right diffusive limit both as $\eps\rightarrow 0$ and $\sigma$ relatively large. Besides, it is better for numerical boundary treatments, e.g., the inflow-outflow close-loop boundary condition for the Marshak wave problem in Section \ref{sec4}.
\end{rem}
%
\section{Numerical methods}
\label{sec3}
\setcounter{equation}{0}
\setcounter{figure}{0}
\setcounter{table}{0}

In this section, we will describe the discontinuous Galerkin (DG) finite element spatial discretization for the system \eqref{eq:mmd2}, where local DG is used for the diffusive terms. We take the two-dimensional case as an example, and the one-dimensional problem can be formulated similarly. In time, the scheme is coupled with high order globally stiffly accurate explicit-implicit (IMEX) methods. The resulting scheme is nonlinear only for the macroscopic variables $\rho$ and $T$. We use a Picard iteration with a predictor-corrector procedure, so that the nonlinear system is decoupled to a linear system, coupled with algebraic nonlinear equations restricted in each element, so that robust and fast convergence can be obtained. 

\subsection{Semi-discrete Discontinuous Galerkin scheme}
\label{sec3.1}
 
We consider a partition of the space computational domain $\mathbb{D}$ with a set of non-overlapping elements $\mathcal{T}_h=\{K\}$, which can completely cover the domain $\mathbb{D}$. $h$ is the maximum edge size of these elements. For simplicity, in this work, we consider a rectangular computational domain $\mathbb{D}$ with a partition of square elements $K=I_{ij}=I_i\times I_j$. Here $I_i=[\xL,\xR]$, $I_j=[\yL,\yR]$ and $(x_i,y_j)$ is the center of each element $K$. 

Given any non-negative integer vector $\bk=(k_1,k_2)$, we define a finite dimensional discrete piecewise polynomial space as follows
\begin{equation}
V_h^\bk=\left\{u\in L^2(\mathbb{D}): \,\, u|_K\in \mathcal{Q}^\bk(K), \,\,\forall \, K\in \mathcal{T}_h\right\},
\label{eq:DiscreteSpace:1mesh}
\end{equation}
where $\mathcal{Q}^\bk(K)$ consists of tensor product polynomials of degree not exceeding $k_i$ along the $i$-th direction on the element $K$, for $i=1,2$. Here $\mathcal{Q}^\bk(K)$ is chosen as it has better $p$-adaptivity along each spatial direction on a rectangular domain. The $\mathcal{P}^k(K)$ local space which consists of polynomials of degree not exceeding $k$ can be used for the following scheme too.

For convenience, on each edge $e$ of $\mathcal{T}_h$, we follow \cite{du2018} to define a unit normal vector $\bn^e$ in the following way. If $e\in\partial \mathbb{D}$, $\bn^e$ is defined as the unit normal vector pointing outside of 
$\mathbb{D}$. For an interior edge $e=\partial K^+\cap\partial K^-$, we first denote $\bn^+$ and $\bn^-$ as the outward unit normal vectors of e taken from the elements $K^+$ and $K^-$, respectively. We fix $\bn^e$ as one of $\bn^\pm$. If we use the notations $u^+$ and $u^-$ to denote the values of a function $u$ on $e$ taken from $K^+$ and $K^-$, respectively, then the
jump notation $[u]$ over an edge $e$ for a scalar valued function $u$ is defined as
\[
[u]|_e = -(u^+\bn^+ + u^-\bn^-)\cdot \bn^e,
\]
and for a vector-valued function $\bq$, the jump $[\bq\cdot\bn]$ is defined as
\[
[\bq\cdot\bn]|_e = -(\bq^+\cdot\bn^+ + \bq^-\cdot\bn^-).
\]
Accordingly, the averages of $u$ and $\bq\cdot\bn$ are defined as
\[
\{u\}|_e = -\frac12(u^+\bn^+ - u^-\bn^-)\cdot \bn^e,\;\;\{\bq\cdot\bn\}|_e = -\frac12(\bq^+\cdot\bn^+ - \bq^-\cdot\bn^-).
\]
Taking square elements as an example, along the $x$ direction, $K^-$ and $K^+$ are the left and right elements associated to the edge $e$ respectively. We may fix $\bn^e=\bn^-$, then
\[
 [u]|_e=u^+-u^-,\;\; [\bq\cdot\bn]|_e=\bq^+-\bq^-, \;\;\{u\}|_e=\frac12(u^++u^-),\;\; \{\bq\cdot\bn\}|_e=\frac12(\bq^++\bq^-).
\]

With these notations, let $\bq=\nabla \rho$, a semi-discrete (local) DG scheme for \eqref{eq:mmd2} is given below. Looking for $\rho_h(t,\vr)$, $g_h(t,\vr,\vO)$, $T_h(t,\vr) \in V_h^{\bk}$ and $\bq_h(t,\vr)=(q_{1,h}(t,\vr),q_{2,h}(t,\vr))^T$ where $q_{i,h}(t,\vr) \in V_h^{\bk}$ for $i=1,2$, such that $\forall \, \testR(\vr)$, $\testG(\vr), \testE(\vr) \in V_h^{\bk}$, and $\testP(\vr)=(z_1(\vr),z_2(\vr))^T$ where $z_i(\vr) \in V_h^{\bk}$ for $i=1,2$, we have
\beq
\label{eq:SDG}
\left\{
\begin{array}{l}
	(\df_t \rho_h, \testR)+\frac{c}{\ssig}\,a_h(\langle \vO\,g_h \rangle, \testR) + c\,\frac{\omega}{3\sigma_0} \,d_h(\bq_h,\testR) =  -\frac{c}{\reps^2}(\rsig \,(\rho_h-\Phi_h),\testR)+ c\,\frac{\omega}{3\sigma_0} \,d_h(\bq_h,\testR), \\ \, \\

	\reps^2\,\tC \,(\df_t T_h, \testE) = (\rsig\,(\rho_h-\Phi_h),\testE), \\	\, \\(\bq_h,\testP) = G_h(\rho_h,\testP), \\ \, \\
	(\df_t g_h, \testG)+ \frac{c/\ssig}{\reps}\, b_h(\vO\,g_h, \testG)+\frac{c/\ssig}{\reps^2} r_h(\vO\,\rho_h, \testG)	
	=-\frac{c}{\reps^2}(\rsig \,g_h, \testG),	
\end{array}
\right.
\eeq
where
\beq
\label{eq:DGh}
\left\{
\begin{array}{l}
	a_h(\langle \vO\,g_h \rangle,\testR)=-\sum_{K\in\mathcal{T}_h} \int_K \langle \vO\,g_h \rangle \divx\testR(\vr) d\vr + \sum_{K\in\mathcal{T}_h} \int_{\partial K} \bn^e\cdot \widehat{\langle \vO\,g_h\rangle}\,\testR(\vr)\,ds,\\ \, \\
	d_h(\bq_h, \testR)=-\sum_{K\in\mathcal{T}_h} \int_K \bq_h \divx\testR(\vr) d\vr + \sum_{K\in\mathcal{T}_h} \int_{\partial K} \bn^e\cdot \widehat{\bq_h}\,\testR(\vr)\,ds,\\ \, \\
	b_h(\vO\,g_h,\testG)=((\mI-\Pi)\mD_h(g_h; \vO), \testG(\vr)) =(\mD_h(g_h; \vO) - \langle\mD_h(g_h; \vO)\rangle, \testG(\vr)),\\ \, \\
	r_h(\vO\,\rho_h, \testG)=-\sum_{K\in\mathcal{T}_h} \int_K \rho_h\,\vO\, \divx\testG(\vr) d\vr + \sum_{K\in\mathcal{T}_h} \int_{\partial K} \bn^e\cdot \vO\,\widehat{\rho_h}\,\testG(\vr)\,ds,\\ \, \\
	G_h(\rho_h, \testP)=-\sum_{K\in\mathcal{T}_h} \int_K \rho_h \nabla\testP(\vr) d\vr + \sum_{K\in\mathcal{T}_h} \int_{\partial K}  \widehat{\rho_h}\,\testP(\vr)\,ds.
\end{array}
\right.
\eeq
Here and below, the standard inner product $(\cdot, \cdot)$ for the $L^2(\mathbb{D})$ space is used, see e.g. in \eqref{eq:SDG}, the first term in the first equation and two terms in the second equation. The function $\mD_h(g_h; \vO)$ in the third equation of \eqref{eq:DGh} is defined as an upwind discretization of $\divx(\vO\,g)$ within the DG framework,
\begin{align}
(\mD_h(g_h; \vO), \testG)
=-\sum_{K\in\mathcal{T}_h} \int_K \vO\,g_h \divx\testG(\vr) d\vr + \sum_{K\in\mathcal{T}_h} \int_{\partial K} \bn^e\cdot \left(\widetilde{\vO\,g_h}\testG(\vr)\right)\,ds,
\label{eq:mD}
\end{align}
with $\widetilde{\vO\,g}$ being an upwind numerical flux consistent to $\vO\,g$,
\beq
\label{eq:vg:upwind:L-1}
\widetilde{\vO\,g}:=
\left\{
\begin{array}{ll}
	\vO\, g^-,&\mbox{if}\; \vO\cdot\bn^e\,>0,\\
	\vO\, g^+,&\mbox{if}\; \vO\cdot\bn^e\,<0.
\end{array}
\right.
\eeq
Remember that $\bn^e$ is a fixed unit normal vector of $\bn^-$ and $\bn^+$ associated to two
elements $K^-$ and $K^+$ and $e=\partial K^+\cap\partial K^-$. $g^\pm$ are the values of $g$ on $e$ taken from $K^\pm$ respectively. $\bn^e\cdot \widehat{\langle \vO\,g_h\rangle}$, $\bn^e\cdot \widehat{\bq_h}$ and $\widehat{\rho_h}$ are also numerical fluxes, and they are consistent to the physical fluxes $\bn^e\cdot\langle \vO\,g\rangle$, $\bn^e\cdot\bq$ and $\rho$, respectively. The following choices could be taken:
\begin{subequations}
	\label{eq:flux}
	\begin{align}
	\label{eq:flux:1}
	{\textrm{alternating left-right:}}  &\;\; \bn^e\cdot \widehat{\langle \vO\,g\rangle} = \bn^-\cdot{\langle \vO\,g\rangle}^-,\; \bn^e \cdot \widehat{\bq}=\bn^- \cdot \widehat{\bq}^-, \; \hat{\rho} = {\rho}^+; \; \\
	\textrm{alternating right-left:} &\;\; \bn^e\cdot \widehat{\langle \vO\,g\rangle} = -\bn^+\cdot{\langle \vO\,g\rangle}^+,\; \bn^e \cdot \widehat{\bq}=-\bn^+ \cdot \widehat{\bq}^+, \; \hat{\rho} = {\rho}^-; \;\\
	\label{eq:flux:2}
	{\textrm{central:}}  &\;\;\bn^e\cdot \widehat{\langle \vO\,g\rangle} = \{\bn\cdot{\langle \vO\,g\rangle}\},\; \bn^e \cdot \widehat{\bq}=\{\bn \cdot \widehat{\bq}\}, \; \hat{\rho} = \{\rho\}\;.
	\end{align}
\end{subequations}

\subsection{Temporal discretization}
\label{sec3.2}

To ensure the correct asymptotic property as $\reps\rightarrow0$, the semi-discrete method in \eqref{eq:SDG} is further coupled with globally stiffly accurate IMEX Runge-Kutta methods in time 
\cite{ascher1997implicit,boscarino2011implicit}. 

We first give a fully discrete scheme with a {\em first order} IMEX method, termed as DG-IMEX1. Given $\rho^n_h(\cdot)$, $g^n_h(\cdot,\vO)$, $T^n_h(\cdot) \in V_h^{\bk}$ and $\bq^n_h(\cdot)=(q^n_{1,h}(\cdot),q^n_{2,h}(\cdot))^T$ where $q^n_{i,h}(\cdot) \in V_h^{\bk}$ for $i=1,2$, which approximate the solution $\rho$, $g$, $T$ and $\bq$ at $t=t^n$, we look for $\rho^{n+1}_h(\cdot)$, $g^{n+1}_h(\cdot,\vO)$, $T^{n+1}_h(\cdot) \in V_h^{\bk}$ and $\bq^{n+1}_h(\cdot)=(q^{n+1}_{1,h}(\cdot),q^{n+1}_{2,h}(\cdot))^T$ where $q^{n+1}_{i,h}(\cdot) \in V_h^{\bk}$, such that
$\forall \, \testR(\vr)$, $\testG(\vr), \testE(\vr) \in V_h^{\bk}$, and $\testP(\vr)=(z_1(\vr),z_2(\vr))^T$ where $z_i(\vr) \in V_h^{\bk}$ for $i=1,2$, we have
\begin{subequations}
	\label{eq:DG}
	\begin{align}
\label{eq:DG:a}	&\left(\frac{\rho^{n+1}_h-\rho^n_h}{\Delta t}, \testR\right)+\frac{c}{\ssig}\,a_h(\langle \vO\,g^n_h \rangle, \testR) + c\,\frac{\omega}{3\sigma_0} \,d_h(\bq^n_h,\testR) = \\ \, \nonumber\\ &\hspace{5cm}-\frac{c}{\reps^2}(\rsig\,(\rho^{n+1}_h-\Phi^{n+1}_h),\testR)+c\,\frac{\omega}{3\sigma_0} \,d_h(\bq^{n+1}_h,\testR), \nonumber\\ \,\nonumber \\
\label{eq:DG:b}	&\reps^2\,\tC \,\left(\frac{T^{n+1}_h-T^n_h}{\Delta t}, \testE\right) = \left(\rsig \,(\rho^{n+1}_h-\Phi^{n+1}_h),\testE\right), \\ \,\nonumber \\
\label{eq:DG:c}	&(\bq^{n+1}_h,\testP) = G_h(\rho^{n+1}_h,\testP),\\	\, \nonumber\\
\label{eq:DG:d}	&\left(\frac{g^{n+1}_h-g^n_h}{\Delta t}, \testG\right)+ \frac{c/\ssig}{\reps}\, b_h(\vO\,g^n_h, \testG)+\frac{c/\ssig}{\reps^2}\, r_h(\vO\,\rho^{n+1}_h, \testG)
	=-\frac{c}{\reps^2}\,(\rsig\,g^{n+1}_h, \testG).	
\end{align}
\end{subequations}
As we can see that the most stiff terms, in both the convective and collisional terms, with a scale of $\frac{1}{\reps^2}$, are treated implicitly. The two added diffusive terms, one is treated explicitly and the other is treated implicitly. 

For the scheme \eqref{eq:DG}, we may first solve the first three equations to update $\rho^{n+1}_h$, $T^{n+1}_h$ and $\bq^{n+1}_h$, then we substitute $\rho^{n+1}_h$ and $T^{n+1}_h$ to the last equation to update $g^{n+1}_h$. Notice that $\rsig$ might depend on $T$ nonlinearly. In this case, the only unknown $g^{n+1}_h$ involving the angular direction $\vO$ is updated explicitly, which is another main point for the design of asymptotic preserving scheme under the micro-macro decompsition framework \cite{Jang2015high}. However, we would note that due to the nonlinear relation
$\Phi=acT^4/|\vO|$, the first three equations of \eqref{eq:DG} still appears to be a fully nonlinear coupled system for the unknowns $\rho^{n+1}_h$ and $T^{n+1}_h$, this is the main difference from the linear radiative transport equations in \cite{Jang2015high,Peng2020stability}. Before we talk about how to deal with this nonlinearity, let us finish our high order temporal discretization.

To achieve high order accuracy in time, we adopt the globally stiffly accurate multi-stage IMEX Runge-Kutta (RK) schemes \cite{boscarino2011implicit,Jang2015high,Peng2020stability}. A multi-stage IMEX RK scheme can be represented by a double Butcher tableau
\begin{equation}\label{eq: B_table}
\begin{array}{c|c}
\tilde{c} & \tilde{\mathcal{A}}\\
\hline
\vspace{-0.25cm}
\\
& \tilde{b^T} \end{array}, \ \ \ \ \  
\begin{array}{c|c}
{c} & {\mathcal{A}}\\
\hline
\vspace{-0.25cm}
\\
& {b^T} \end{array},
\end{equation}
where both $\tilde{\mathcal{A}} = (\tilde{a}_{ij})$ and $\mathcal{A} = (a_{ij})$ are $(s+1)\times (s+1)$ matrices, with $\tilde{\mathcal{A}}$ being lower triangular with zero diagonal entries. The coefficients $\tilde{c}$ and $c$ are given by the usual relation $\tilde{c}_i = \sum_{j=0}^{i-1} \tilde a_{ij}$, $c_i = \sum_{j=0}^{i} a_{ij}$,
and vectors $\tilde{b} = (\tilde{b}_j)$ and $b = (b_j)$ provide the quadrature weights to combine internal stages of the RK method. The IMEX RK scheme is said to be {\em globally stiffly accurate} \cite{boscarino2011implicit} if
\begin{equation}
c_s = \tilde{c}_s = 1, \;\;\textrm{and}\;\;  a_{sj} = b_j, \;\; \tilde{a}_{sj} = \tilde{b}_j, \;\;\forall j=0, \cdots, s.
\label{eq:gsa}
\end{equation}
This double Butcher tableau corresponds to explicit and implicit time discretization respectively. In this work, the type ARS IMEX RK scheme is used, where $\mathcal{A}$ in \eqref{eq: B_table} has the following structure,
\[
\left[
	\begin{array}{cc}
	0 & 0 \\
	0 & \hat{\mathcal{A}}
	\end{array}
\right]
\]
where the diagonal entries of $\hat{\mathcal{A}}$ are nonzero. 

Now we will apply a general globally stiffly accurate IMEX RK method, represented by \eqref{eq: B_table} with the property \eqref{eq:gsa}, to the semi-discrete DG scheme \eqref{eq:SDG}. For a $p$-th order IMEX RK method, our scheme is termed as ``DG-IMEXp". This is the same implicit-explicit strategy as in the first order case. Given $\rho^n_h(\cdot)$, $g^n_h(\cdot,\vO)$, $T^n_h(\cdot) \in V_h^{\bk}$ and $\bq^n_h(\cdot)=(q^n_{1,h}(\cdot),q^n_{2,h}(\cdot))^T$ where $q^n_{i,h}(\cdot) \in V_h^{\bk}$ for $i=1,2$ that approximate the solution $\rho$, $g$, $T$ and $\bq$ at $t=t^n$, we look for $\rho^{n+1}_h(\cdot)$, $g^{n+1}_h(\cdot,\vO)$, $T^{n+1}_h(\cdot) \in V_h^{\bk}$ and $\bq^{n+1}_h(\cdot)=(q^{n+1}_{1,h}(\cdot),q^{n+1}_{2,h}(\cdot))^T$ where $q^{n+1}_{i,h}(\cdot) \in V_h^{\bk}$, such that
$\forall \, \testR(\vr)$, $\testG(\vr), \testE(\vr) \in V_h^{\bk}$, and $\testP(\vr)=(z_1(\vr),z_2(\vr))^T$ where $z_i(\vr) \in V_h^{\bk}$ for $i=1,2$, we have
\begin{subequations}
\label{eq:FDG:GSA}
\begin{align}
	\label{eq:FDG:GSA:a}&\left(\frac{\rho^{n+1}_h-\rho^n_h}{\Delta t}, \testR\right)+\frac{c}{\ssig}\,\sum_{\ell=0}^s \tilde{b}_\ell \,a_h(\langle \vO\,g^{(\ell)}_h \rangle, \testR) + c\,\frac{\omega}{3\sigma_0} \,\sum_{\ell=0}^s \tilde{b}_\ell \,d_h(\bq^{(\ell)}_h,\testR) = \\  & \hspace{5cm} -\frac{c}{\reps^2}\sum_{\ell=0}^s b_\ell \, (\rsig\,(\rho^{(\ell)}_h-\Phi^{(\ell)}_h),\testR)+c\,\frac{\omega}{3\sigma_0} \,\sum_{\ell=0}^s \,b_\ell\,d_h(\bq^{(\ell)}_h,\testR), \nonumber\\ \, \nonumber\\
	\label{eq:FDG:GSA:b}&\reps^2\,\tC \,\left(\frac{T^{n+1}_h-T^n_h}{\Delta t}, \testE\right) = \sum_{\ell=0}^s \,b_\ell\,\left(\rsig \,(\rho^{(\ell)}_h-\Phi^{(\ell)}_h),\testE\right), \\ \, \nonumber\\
	\label{eq:FDG:GSA:c}&(\bq^{n+1}_h,\testP) = G_h(\rho^{n+1}_h,\testP),\\	\, \nonumber\\
	\label{eq:FDG:GSA:d}&\left(\frac{g^{n+1}_h-g^n_h}{\Delta t}, \testG\right)+ \frac{c/\ssig}{\reps} \,\sum_{\ell=0}^s \tilde{b}_\ell \,b_h(\vO\,g^{(\ell)}_h, \testG)+\frac{c/\ssig}{\reps^2} \,\sum_{\ell=0}^s b_\ell \,r_h(\vO\,\rho^{(\ell)}_h, \testG)
	=\\  &\hspace{5cm}-\frac{c/\ssig}{\reps^2}\,\sum_{\ell=0}^s b_\ell \,(\rsig\,g^{(\ell)}_h, \testG). 	\nonumber
\end{align}
\end{subequations}
The approximations at the internal stages $\rho_h^{(\ell)}(\cdot)$, $g_h^{(\ell)}(\cdot,\vO)$, $T_h^{(\ell)}(\cdot) \in V_h^\bk$ and $\bq_h^{(\ell)}(\cdot)=(q_{1,h}^{(\ell)}(\cdot),q_{2,h}^{(\ell)}(\cdot))$ where $q_{i,h}^{(\ell)}(\cdot)\in V^\bk_h$ for $i=1,2$, with $\ell=1, \cdots, s$, would satisfy
\begin{subequations}
	\label{eq:FDG:GSA:stage}
	\begin{align}
	\label{eq:FDG:GSA:stage:a}&\left(\frac{\rho^{(\ell)}_h-\rho^n_h}{\Delta t}, \testR\right)+\frac{c}{\ssig}\,\sum_{j=0}^{\ell-1} \tilde{a}_{\ell,j} \,a_h(\langle \vO\,g^{(j)}_h \rangle, \testR) + c\,\frac{\omega}{3\sigma_0} \,\sum_{j=0}^{\ell-1} \tilde{a}_{\ell,j} \,d_h(\bq^{(j)}_h,\testR) = \\ &\hspace{5cm} -\frac{c}{\reps^2}\sum_{j=0}^{\ell} a_{\ell,j} \, (\rsig \, (\rho^{(j)}_h-\Phi^{(j)}_h),\testR)+c\,\frac{\omega}{3\sigma_0} \,\sum_{j=0}^{\ell} \,a_{\ell,j}\,d_h(\bq^{(j)}_h,\testR), \nonumber\\ \, \nonumber\\
	\label{eq:FDG:GSA:stage:b}&\reps^2\,\tC \,\left(\frac{T^{(\ell)}_h-T^n_h}{\Delta t}, \testE\right) = \sum_{j=0}^{\ell} \,a_{\ell,j}\,\left(\rsig \,(\rho^{(j)}_h-\Phi^{(j)}_h),\testE\right), \\ \, \nonumber\\
	\label{eq:FDG:GSA:stage:c}&(\bq^{(\ell)}_h,\testP) = G_h(\rho^{(\ell)}_h,\testP),\\	\, \nonumber\\
	\label{eq:FDG:GSA:stage:d}&\left(\frac{g^{(\ell)}_h-g^n_h}{\Delta t}, \testG\right)+ \frac{c/\ssig}{\reps} \,\sum_{j=0}^{\ell-1} \tilde{a}_{\ell,j} \,b_h(\vO\,g^{(j)}_h, \testG)+\frac{c/\ssig}{\reps^2} \,\sum_{j=0}^{\ell} a_{\ell,j} \,r_h(\vO\,\rho^{(j)}_h, \testG)
	= \\  &\hspace{5cm}-\frac{c}{\reps^2}\,\sum_{j=0}^{\ell} a_{\ell,j} \,(\rsig\,g^{(j)}_h, \testG). \nonumber	
	\end{align}
\end{subequations}
Notice that for $\ell=0$, $\rho^{(0)}_h=\rho^n_h,\;g^{(0)}=g^n_h,\;T^{(0)}_h=T^n_h$ and $\bq^{(0)}_h=\bq^n_h$. 
This high order temporal discretization \eqref{eq:FDG:GSA}-\eqref{eq:FDG:GSA:stage} can be solved similarly as the first order scheme \eqref{eq:DG} in a stage-by-stage manner. That is,
in each internal stage, from \eqref{eq:FDG:GSA:stage}, we first solve the first three equations to update $\rho_h^{(\ell)}(\cdot)$, $T_h^{(\ell)}(\cdot)$ and $\bq_h^{(\ell)}(\cdot)$, and then substituting $\rho_h^{(\ell)}(\cdot)$ and $T_h^{(\ell)}(\cdot)$ into the fourth equation to get $g_h^{(\ell)}(\cdot,\vO)$. 

After all stage values are obtained, the numerical solutions $\rho^{n+1}_h(\cdot)$, $g^{n+1}_h(\cdot,\vO)$, $T^{n+1}_h(\cdot) $ and $\bq^{n+1}_h(\cdot)$ at the next time level $t^{n+1}$ can be accumulated from \eqref{eq:FDG:GSA}. For a globally stiffly accurate IMEX RK scheme, from \eqref{eq:gsa}, the solutions at $t^{n+1}$ concide with the last stage values, that is, 
\beq
\rho^{n+1}_h=\rho^{(s)}_h, \quad g^{n+1}_h=g^{(s)}_h, \quad T^{n+1}_h=T^{(s)}_h \; \textrm{ and } \; \Phi^{n+1}=ac(T^{n+1}_h)^4/|\vO|.
\eeq
However, numerically we find that due to \eqref{eq:FDG:GSA:stage:a}-\eqref{eq:FDG:GSA:stage:c} is solved iteratively (as described in the next section), the exact conservation is lost and we cannot get the right front propagation of the material temperature. Keeping $g^{n+1}_h=g^{(s)}_h$, combining \eqref{eq:FDG:GSA:a} and \eqref{eq:FDG:GSA:b} we propose to update $\rho^{n+1}_h$ and $T^{n+1}_h$ again with
\begin{subequations}
\label{eq:FDG:GSA:n+1}
\begin{align}
&\left(\frac{\rho^{n+1}_h-\rho^n_h}{\Delta t}, \testR\right)+\tC \,\left(\frac{T^{n+1}_h-T^n_h}{\Delta t} , \testR\right) = \\ & \hspace{2cm} +\frac{c}{\ssig}\,\sum_{\ell=0}^s \tilde{b}_\ell \,a_h(\langle \vO\,g^{(\ell)}_h \rangle, \testR) + c\,\frac{\omega}{3\sigma_0} \,\sum_{\ell=0}^s \tilde{b}_\ell \,d_h(\bq^{(\ell)}_h,\testR)+c\,\frac{\omega}{3\sigma_0} \,\sum_{\ell=0}^s \,b_\ell\,d_h(\bq^{(\ell)}_h,\testR), \nonumber \\ \, \nonumber \\
&\reps^2\,\tC \,\left(\frac{T^{n+1}_h-T^n_h}{\Delta t}, \testE\right) - b_s\,\left(\rsig \,(\rho^{n+1}_h-\Phi^{n+1}_h),\testE\right) = \sum_{\ell=0}^{s-1} \,b_\ell\,\left(\rsig \,(\rho^{(\ell)}_h-\Phi^{(\ell)}_h),\testE\right). 
\end{align}
\end{subequations}
The idea is that spatial derivatives are discretized in a conservative flux difference form for updating $\rho^{n+1}_h$. With \eqref{eq:FDG:GSA:n+1}, both conservation and $\rho^{n+1}_h=\Phi^{n+1}_h+\mO(\reps)$ when $\reps\ll 1$ can be ensured.

\subsection{Picard iteration with a predictor-corrector procedure}
\label{sec3.3}
As mentioned in the last subsection, the first three equations in the scheme \eqref{eq:DG} or \eqref{eq:FDG:GSA:stage} are fully coupled nonlinear system. To avoid solving the nonlinear system directly, we empoly a Picard iteration with a predictor-corrector procedure. This is motivated from the work in \cite{tang2020}, where this procedure is designed for an accurate front capturing asymptotic preserving scheme. The idea is that, to avoid a coupled global nonlinear system, we need to seperate the spatial differential operator $\Delta \rho$ from the nonlinear term $\Phi=acT^4/|\vO|$ with respect to $T$. 

We start with the {\em first order} in time scheme \eqref{eq:DG}. We focus on the first three equations and we rewrite them in the following way
\begin{subequations}
	\label{eq:DG:picard}
	\begin{align}
	\label{eq:DG:picard:a}&(\rho^{n+1}_h, \testR)+c\,\tC \,\left(T^{n+1}_h, \testR\right) -c\,\Delta t\,\frac{\omega}{3\sigma_0} \,d_h(\bq^{n+1}_h,\testR)\;=\;(\textrm{RHS}_1^n,\testR), \\ \,\nonumber \\
	\label{eq:DG:picard:b}&\reps^2\,\tC \,\left(T^{n+1}_h, \testE\right) -\Delta t\, \left(\rsig\,(\rho^{n+1}_h-\Phi^{n+1}_h),\testE\right)\;=\;(\textrm{RHS}_2^n,\testE), \\ \,\nonumber \\
	\label{eq:DG:picard:c}&(\bq^{n+1}_h,\testP) - G_h(\rho^{n+1}_h,\testP)\;=\; 0.	
	\end{align}
\end{subequations}
We have combined \eqref{eq:DG:a} and \eqref{eq:DG:b} to obtain \eqref{eq:DG:picard:a}. 
Besides, we have put all unknown values at the time level $t^{n+1}$ on the left, while all values at the time level $t^n$ on the right, where 
\[
(\textrm{RHS}_1^n,\testR)=(\rho^{n}_h, \testR)+c\,\tC \,\left(T^{n}_h, \testR\right) -c\,\Delta t\,\frac{\omega}{3\sigma_0} \,d_h(\bq^{n}_h,\testR)-\frac{c}{\ssig}\,\Delta t\,a_h(\langle\vO\,g^n_h\rangle,\testR),
\]
and
\[
(\textrm{RHS}_2^n,\testE)=\reps^2\,\tC \,\left(T^{n}_h, \testE\right).
\]
Now the Picard iteration for solving \eqref{eq:DG:picard} is defined as follows:
for the iterative number $m$ starting at $m=0$, where $\rho^{n+1,(0)}_h=\rho^n_h$,
$T^{n+1,(0)}_h=T^n_h$ and $\bq^{n+1,(0)}=\bq^n_h$, we update the unknowns $\rho^{n+1,(m+1)}_h$,
$T^{n+1,(m+1)}_h$ and $\bq^{n+1,(m+1)}$ iteratively by the following two steps:
\begin{itemize}
\item Step 1: we first solve the following linear system to obtain $\rho^{n+1,(m+1)}_h$ and $\bq^{n+1,(m+1)}$
\begin{equation}
\label{eq:step1}
\left\{
\begin{array}{l}
	(\rho^{n+1,(m+1)}_h, \testR) -c\,\Delta t\,\frac{\omega}{3\sigma_0} \,d_h(\bq^{n+1,(m+1)}_h,\testR)\;=\;(\textrm{RHS}_1^n,\testR)-c\,\tC \,\left(T^{n+1,(m)}_h, \testR\right), \\ \,\\
	(\bq^{n+1,(m+1)}_h,\testP) - G_h(\rho^{n+1,(m+1)}_h,\testP)\;=\; 0.	
\end{array}
\right.
\end{equation}
Notice that $d_h(\bq^{n+1,(m+1)}_h,\testR)$ combined with $(\bq^{n+1,(m+1)}_h,\testP)$ gives the local DG discretization for the Laplacian $\Delta \rho$, which results in a linear system;
\item Step 2: with $\rho^{n+1,(m+1)}_h$ and $\bq^{n+1,(m+1)}$ obtained from Step 1, we now solve
\beq
\label{eq:step2}
\left\{
\begin{array}{l}
	(\rho^{n+1,(m+1)}_h, \testR)+c\,\tC \,\left(T^{n+1,(m+1)}_h, \testR\right) \;=\;(\textrm{RHS}_1^n,\testR)+c\,\Delta t\,\frac{\omega}{3\sigma_0} \,d_h(\bq^{n+1,(m+1)}_h,\testR), \\ \, \\
	\reps^2\,\tC \,\left(T^{n+1,(m+1)}_h, \testE\right) -\Delta t\, \left(\rsig\,(\rho^{n+1,(m+1)}_h-\Phi^{n+1,(m+1)}_h),\testE\right)\;=\;(\textrm{RHS}_2^n,\testE).	
\end{array}
\right.
\eeq
$d_h(\bq^{n+1,(m+1)}_h,\testR)$ can be precomputed from $\bq^{n+1,(m+1)}_h$. Instead of using $\rho^{n+1,(m+1)}_h$ obtained from step 1, here we express it by $T_h^{n+1,(m+1)}$ and $\Phi_h^{n+1,(m+1)}$ from the second equation of \eqref{eq:step2}, and substitute it into the first equation. In this way, the system \eqref{eq:step2} becomes a local algebraic nonlinear system for $T_h^{n+1,(m+1)}$ within each element $K$, so that solving a global nonlinear system is avoid. The Newton iteration can be simply used to solve this algebraic nonlinear system.
\end{itemize}
In Step 2, if we assume $\rsig$ is simply a positive constant, from the second equation, we have
\beq
(\rho^{n+1,(m+1)}_h,\testE)\;=\;(\Phi^{n+1,(m+1)}_h,\testE)-\reps^2\,\tC \,\left(T^{n+1,(m+1)}_h, \testE\right)/(\rsig\,\Delta t)-(\textrm{RHS}_2^n,\testE)/(\rsig\,\Delta t).
\eeq
substituting it into the first equation, it yields
\beq
\label{eq:T}
(\Phi^{n+1,(m+1)}_h,\testR)+c\,\tC \,\left(T^{n+1,(m+1)}_h, \testR\right)-\reps^2\,\tC \,\left(T^{n+1,(m+1)}_h, \testR\right)/(\rsig\,\Delta t) \;=\;(\textrm{RHS}_3^n,\testR),
\eeq
where
\[
(\textrm{RHS}_3^n,\testR)=(\textrm{RHS}_1^n,\testR)+c\,\Delta t\,\frac{\omega}{3\sigma_0} \,d_h(\bq^{n+1,(m+1)}_h,\testR)+(\textrm{RHS}_2^n,\testR)/(\rsig\,\Delta t).
\] 
\eqref{eq:T} in the diffusive limit as $\reps\ra0$, it turns out to be 
\beq
\label{eq:TD}
\begin{array}{l}
(\Phi^{n+1,(m+1)}_h,\testR)+c\,\tC \,\left(T^{n+1,(m+1)}_h, \testR\right) \;=\;(\rho^{n}_h, \testR)+c\,\tC \,\left(T^{n}_h, \testR\right) -c\,\Delta t\,\frac{\omega}{3\sigma_0} \,d_h(\bq^{n}_h,\testR) \\ \hspace{6.5cm} -\frac{c}{\ssig}\,\Delta t\,a_h(\langle\vO\,g^n_h\rangle,\testR)+c\,\Delta t\,\frac{\omega}{3\sigma_0} \,d_h(\bq^{n+1,(m+1)}_h,\testR).
\end{array}
\eeq
Since $\rho^n_h$ approaches $\Phi^n_h$ from \eqref{eq:DG:picard:b}, combining with \eqref{eq:DG:d} in the limit $\reps\ra0$, \eqref{eq:TD} is
a consistent discretization for the diffusive equation \eqref{diffeq} (for more details, see the asymptotic analysis in the following section). Namely \eqref{eq:TD} is mimicking a numerical discretization for the diffusive limiting equation \eqref{diffeq}. 

The two steps are solved with the iterative number $m$ untill convergent, where the stop criteria is defined as
\[
\|\rho^{n+1,(m+1)}_h-\rho^{n+1,(m)}_h\|_2 < \delta.
\]
The $L^2$ norm is used and we take $\delta=10^{-9}$ in our numerical tests.

This procedure can be similarly applied to high order IMEX RK time discretizations \eqref{eq:FDG:GSA:stage} at each stage. That is, for the first three equations in \eqref{eq:FDG:GSA:stage}, at the stage $\ell$ ($\ell=1,\cdots,s$), we rewrite them as
\begin{subequations}
	\label{eq:FDG:GSA2}
	\begin{align}
	\label{eq:FDG:GSA2:a}&(\rho^{(\ell)}_h, \testR)+c\,\tC \,\left(T^{(\ell)}_h, \testR\right) -c\,a_{\ell\ell}\,\Delta t\,\frac{\omega}{3\sigma_0} \,d_h(\bq^{(\ell)}_h,\testR)\;=\;(\textrm{RHS}_1,\testR), \\ \,\nonumber \\
	\label{eq:FDG:GSA2:b}&\reps^2\,\tC \,\left(T^{(\ell)}_h, \testE\right) -a_{\ell\ell}\,\Delta t\, \left(\rsig \,(\rho^{(\ell)}_h-\Phi^{(\ell)}_h),\testE\right)\;=\;(\textrm{RHS}_2,\testE), \\ \,\nonumber \\
	\label{eq:FDG:GSA2:c}&(\bq^{(\ell)}_h,\testP) - G_h(\rho^{(\ell)}_h,\testP)\;=\; 0.	
	\end{align}
\end{subequations}
Similarly we have combined \eqref{eq:FDG:GSA:stage:a} and \eqref{eq:FDG:GSA:stage:b} to obtain \eqref{eq:FDG:GSA2:a}. All known values before stage $\ell$ are put on the right, where 
\[
(\textrm{RHS}_1,\testR)=(\rho^{n}_h, \testR)+c\,\tC \,\left(T^{n}_h, \testR\right) -\frac{c}{\ssig}\,\Delta t\,\sum_{j=1}^{\ell-1} \tilde{a}_{\ell,j} \,a_h(\langle \vO\,g^{(j)}_h \rangle, \testR) - c\,\Delta t\,\frac{\omega}{3\sigma_0} \,\sum_{j=0}^{\ell-1} \left(\tilde{a}_{\ell,j}-a_{\ell,j}\right) \,d_h(\bq^{(j)}_h,\testR),
\]
and
\[
(\textrm{RHS}_2,\testE)=\reps^2\,\tC \,\left(T^{n}_h, \testE\right)+\Delta t\,\sum_{j=0}^{\ell-1} \,a_{\ell,j}\,\left(\rsig \,(\rho^{(j)}_h-\Phi^{(j)}_h),\testE\right).
\]
Notice that \eqref{eq:FDG:GSA2} is in the same form as \eqref{eq:DG:picard}, so that it can be solved similarly with the Picard iteration described above. The solutions $\rho^{n+1}_h$ and $T^{n+1}_h$ are updated by \eqref{eq:FDG:GSA:n+1} following the same procedure in Step 2. We omit the details here to save space. 


\section{Formal asymptotic preserving analysis}
\label{sec3.5}
In this section, we will perform formal asymptotic preserving (AP) and asymptotically accurate (AA) analyses for the proposed schemes by assuming $\reps \ll 1$, while the mesh parameters $h$ and $\Delta t$ are fixed. We will show that our DG-IMEX1 scheme \eqref{eq:DG} is AP in the sence that when $\reps\ra0$, it becomes a consistent discretization for the limiting diffusive equation \eqref{diffeq}. On the other hand, our DG-IMEXp scheme \eqref{eq:FDG:GSA}-\eqref{eq:FDG:GSA:stage} is asymptotically accurate, namely, the schemes of \eqref{eq:FDG:GSA}-\eqref{eq:FDG:GSA:stage} as $\reps\ra 0$, maintains its order of temporal accuracy for the limiting equation \eqref{diffeq}. The analyses mainly follow the notations and assumptions as used in \cite{Peng2020asymptotic}.

We first assume the initial data is well-prepared without considering the initial layers, namely, $I(t,\vr,\vO)=\rho(t,\vr)+\reps\,g(t,\vr,\vO)$ and $\rho(t,\vr)=\mathcal{O}(1), g(t,\vr,\vO)=\mathcal{O}(1)$. We further make the assumptions that the spatial derivatives of $\rho$ and $g$ at $t=0$, have
comparable scales as $\|\rho(0,\cdot)\|$ and $\|g(0,\cdot,\vO)\|$ with respect to $\reps$, that is, they all are $\mathcal{O}(1)$. Under theses assumptions, all spatial derivative approximations in \eqref{eq:DGh} also all are $\mathcal{O}(1)$. For the small parameters $\reps$, $h$ and $\Delta t$, it is assumed $\reps^2 \ll \Delta t$, $\reps\le\Delta t\le 1$ and $\Delta t/h=\mathcal{O}(1)$ to avoid the explicit dependence on $\Delta t$ and $h$ of the hidden constant in the $\mathcal{O}$ notation.

We begin to write our scheme in a strong form by defining the following linear operators for \eqref{eq:DGh}
\beq
\label{eq:DG:strong1}
(\mD^g_h(\langle \vO\,g_h \rangle), \testR)=a_h(\langle \vO\,g_h \rangle,\testR),\qquad
(\mD^\rho_h(\vO\,\rho_h), \testG)=-r_h(\vO\,\rho_h, \testG),
\eeq
and
\beq
\label{eq:DG:strong2}
(\mD^{grad}_h\rho_h, \testG)=G_h(\rho_h, \testG), \qquad (\mD^{div}_h\bq_h, \testR)=d_h(\bq_h, \testR).
\eeq
They are well-defined bounded operators following the Riesz representation, and determined entirely by the discrete space $V_h^\bk$ and the involved numerical fluxes.

\subsection{Asymptotic preserving (AP) analysis for DG-IMEX1}

The DG-IMEX1 scheme \eqref{eq:DG} in the strong form becomes
\begin{subequations}
	\label{eq:DG:s}
	\begin{align}
	\label{eq:DG:sa}	&\frac{\rho^{n+1}_h-\rho^n_h}{\Delta t}+\frac{c}{\ssig}\,\mD^g_h(\langle \vO\,g^n_h \rangle) + c\,\frac{\omega}{3\sigma_0} \,\mD^{div}_h(\bq^n_h) = -\frac{c}{\reps^2}\rsig\,(\rho^{n+1}_h-\Phi^{n+1}_h)+c\,\frac{\omega}{3\sigma_0} \,\mD^{div}_h(\bq^{n+1}_h), \\ \,\nonumber \\
	\label{eq:DG:sb}	&\reps^2\,\tC \,\frac{T^{n+1}_h-T^n_h}{\Delta t}  = \rsig \,(\rho^{n+1}_h-\Phi^{n+1}_h), \\ \,\nonumber \\
	\label{eq:DG:sc}	&\bq^{n+1}_h = \mD^{grad}_h(\rho^{n+1}_h),\\	\, \nonumber\\
	\label{eq:DG:sd}	&\frac{g^{n+1}_h-g^n_h}{\Delta t}+ \frac{c/\ssig}{\reps}\, (\mI-\Pi)\mD_h(g^n_h;\,\vO)-\frac{c/\ssig}{\reps^2}\, \mD^\rho_h(\vO\,\rho^{n+1}_h)
	=-\frac{c}{\reps^2}\,\rsig\,g^{n+1}_h.	
	\end{align}
\end{subequations}
We notice that $\mD^{div}_h\mD^{grad}_h$ is an approximation of the Laplacian operator $\Delta$. Combining \eqref{eq:DG:sa} and \eqref{eq:DG:sb}, we have
\beq
\label{eq:DG:diff}
\frac{\rho^{n+1}_h-\rho^n_h}{\Delta t}+\frac{c}{\ssig}\,\mD^g_h(\langle \vO\,g^n_h \rangle) + c\,\frac{\omega}{3\sigma_0} \,\mD^{div}_h(\bq^n_h) = -\tC \,\frac{T^{n+1}_h-T^n_h}{\Delta t}+c\,\frac{\omega}{3\sigma_0} \,\mD^{div}_h(\bq^{n+1}_h).
\eeq
We start from the exact solutions at time level $t^n$ and we assume a formal $\reps$-expansion of the quantities $(\rho_h, g_h, \bq_h, T_h, \Phi_h)$, that is
\begin{subequations}
	\label{epsexpansion}
	\begin{align}
	&\rho_h = \rho_{0,h}+\reps\,\rho_{1,h}+\cdots, \\
	&g_h = g_{0,h}+\reps\,g_{1,h}+\cdots, \\
	&\bq_h = \bq_{0,h}+\reps\,\bq_{1,h}+\cdots, \\
	&T_h = T_{0,h}+\reps\,T_{1,h}+\cdots, \\
	&\Phi_h = \Phi_{0,h}+\reps\,\Phi_{1,h}+\cdots, 
	\end{align}
\end{subequations}
Substituting them into \eqref{eq:DG:s}, first from \eqref{eq:DG:sb}, we can see that $\rho_{0,h}^{n+1}=\Phi^{n+1}_{0,h}+\mO(\reps)$, and from \eqref{eq:DG:sd}, we have $g^{n+1}_{0,h}=-\frac{\ssig}{\sigma}\mD^\rho_h(\vO\rho^{n+1}_{0,h})+\mO({\reps})$.
Substituting them with \eqref{eq:DG:sc} into \eqref{eq:DG:sa}, we obtain
\begin{align}
\label{eq:DG:diff2}
&\frac{\Phi^{n+1}_{0,h}-\rho^n_{0,h}}{\Delta t}+c\,\tC \,\frac{T^{n+1}_{0,h}-T^n_{0,h}}{\Delta t} =c\,\mD^g_h(\langle \vO\,\frac{1}{\sigma}\mD^\rho_h(\vO\,\rho^n_{0,h}) \rangle) \\ &\hspace{5.7cm} +c\,\frac{\omega}{3\sigma_0} \,\mD^{div}_h\mD^{grad}_h(\rho^{n+1}_{0,h}-\rho^{n}_{0,h})+\mO(\reps). \nonumber
\end{align}
With the fluxes defined in \eqref{eq:flux}, from \eqref{eq:DGh} we have
\beq
\mD^g_h(\langle \vO\,g^n_h \rangle)=\mD^{div}(\langle \vO\,g^n_h \rangle), \quad \mD^\rho_h(\vO\,\rho^n_h) = \vO\cdot\mD^{grad}(\rho^n_h).
\eeq
Since $\langle\vO\vO\rangle=\frac13\mI$, here $\mI$ is the $2\times 2$ identity matrix, we find that $\mD^g_h(\langle \vO\,\frac{1}{\sigma}\mD^\rho_h(\vO\,\rho^n_{0,h}) \rangle)$ is a local DG approximation to $\nabla\cdot\frac{1}{3\sigma}\nabla\rho$ at time level $t^n$, that is
\beq
\label{eq:diffrho}
\mD^g_h(\langle \vO\,\frac{1}{\sigma}\mD^\rho_h(\vO\,\rho^n_{0,h}) \rangle)=\mD^{div}_h\left(\frac{1}{3\sigma}\mD^{grad}_h\left(\rho^n_{0,h}\right)\right).
\eeq
Due to well-prepared initial conditions, we also have $\rho^n_{0,h}=\Phi^n_{0,h}+\mO(\reps)$. Combining \eqref{eq:DG:diff2} and \eqref{eq:diffrho}, and replacing $\Phi_{0,h}$ by $acT_{0,h}^4/|\vO|$, we obtain
\begin{align}
\label{eq:DG:diff3}
&\frac{a(T^{n+1}_{0,h})^4-a(T^n_{0,h})^4}{\Delta t}+ C_v \,\frac{T^{n+1}_{0,h}-T^n_{0,h}}{\Delta t} =\mD^{div}_h\left(\frac{1}{3\sigma}\mD^{grad}_h\left(ac(T^n_{0,h})^4\right)\right) \\ &\hspace{5cm} +\frac{\omega}{3\sigma_0} \,\mD^{div}_h\mD^{grad}_h(ac(T^{n+1}_{0,h})^4-ac(T^{n}_{0,h})^4)+\mO(\reps), \nonumber
\end{align}
where both sides have been divided by $c/|\vO|$ and $\tC$ is replaced by $C_v$. Since the second term on the right hand side is of order $\mO(\Delta t)$, we can clearly see \eqref{eq:DG:diff3} in the leading order of the unknowns is a consistent discretization for the diffusive limiting equation \eqref{diffeq}, so that the AP property is guaranteed.
\begin{rem}
	We should notice that in the limit scheme \eqref{eq:DG:diff3}, on the right hand side, the second term appears as a linear pernalization
	to the first nonlinear diffusive term, the stability condition is $\Delta t=\mO(h)$ instead of $\Delta t=\mO(h^2)$. For more discussions, see \cite{wang2020}.
\end{rem}

\subsection{Asymptotically accurate (AA) analysis for DG-IMEXp.}
Now let us prove that our DG-IMEXp scheme, as $\reps\ra 0$, is a consistent and high order discretization of the limiting diffusive equation \eqref{diffeq}, especially, the order $p$ of temporal accuracy is maintained without reduction. Similarly, we write the high order DG-IMEXp scheme in a strong form,
\begin{subequations}
	\label{eq:FDG:GSA:s}
	\begin{align}
	\label{eq:FDG:GSA:sa}&\frac{\rho^{n+1}_h-\rho^n_h}{\Delta t}+\frac{c}{\ssig}\,\sum_{\ell=0}^s \tilde{b}_\ell \,\mD^g_h(\langle \vO\,g^{(\ell)}_h \rangle) + c\,\frac{\omega}{3\sigma_0} \,\sum_{\ell=0}^s \tilde{b}_\ell \,\mD^{div}_h(\bq^{(\ell)}_h) = \\  & \hspace{5cm} -\frac{c}{\reps^2}\sum_{\ell=0}^s b_\ell \, \rsig\,(\rho^{(\ell)}_h-\Phi^{(\ell)}_h)+c\,\frac{\omega}{3\sigma_0} \,\sum_{\ell=0}^s \,b_\ell\,\mD^{div}_h(\bq^{(\ell)}_h), \nonumber\\ \, \nonumber\\
	\label{eq:FDG:GSA:sb}&\reps^2\,\tC \,\frac{T^{n+1}_h-T^n_h}{\Delta t} = \sum_{\ell=0}^s \,b_\ell\,\rsig \,(\rho^{(\ell)}_h-\Phi^{(\ell)}_h), \\ \, \nonumber\\
	\label{eq:FDG:GSA:sc}&\bq^{n+1}_h = \mD^{grad}_h(\rho^{n+1}_h),\\	\, \nonumber\\
	\label{eq:FDG:GSA:sd}&\frac{g^{n+1}_h-g^n_h}{\Delta t}+ \frac{c/\ssig}{\reps} \,\sum_{\ell=0}^s \tilde{b}_\ell \,(\mI-\Pi)\mD_h(g^{(\ell)}_h;\,\vO) -\frac{c/\ssig}{\reps^2} \,\sum_{\ell=0}^s b_\ell \,\mD^\rho_h(\vO\,\rho^{(\ell)}_h) 
	= \\ & \hspace{5cm} -\frac{c/\ssig}{\reps^2}\,\sum_{\ell=0}^s b_\ell \,\rsig\,g^{(\ell)}_h, 	\nonumber
	\end{align}
\end{subequations}
with internal stages satisfying
\begin{subequations}
	\label{eq:FDG:GSA:stage:s}
	\begin{align}
	\label{eq:FDG:GSA:stage:sa}&\frac{\rho^{(\ell)}_h-\rho^n_h}{\Delta t} +\frac{c}{\ssig}\,\sum_{j=0}^{\ell-1} \tilde{a}_{\ell,j} \,\mD^g_h(\langle \vO\,g^{(j)}_h \rangle) + c\,\frac{\omega}{3\sigma_0} \,\sum_{j=0}^{\ell-1} \tilde{a}_{\ell,j} \,\mD^{div}_h(\bq^{(j)}_h) = \\  &\hspace{5cm} -\frac{c}{\reps^2}\sum_{j=0}^{\ell} a_{\ell,j} \, (\rsig \, (\rho^{(j)}_h-\Phi^{(j)}_h))+c\,\frac{\omega}{3\sigma_0} \,\sum_{j=0}^{\ell} \,a_{\ell,j}\,\mD^{div}_h(\bq^{(j)}_h), \nonumber\\ \, \nonumber\\
	\label{eq:FDG:GSA:stage:sb}&\reps^2\,\tC \,\frac{T^{(\ell)}_h-T^n_h}{\Delta t} = \sum_{j=0}^{\ell} \,a_{\ell,j}\,\rsig \,(\rho^{(j)}_h-\Phi^{(j)}_h), \\ \, \nonumber\\
	\label{eq:FDG:GSA:stage:sc}&\bq^{(\ell)}_h = \mD^{grad}_h(\rho^{(\ell)}_h),\\	\, \nonumber\\
	\label{eq:FDG:GSA:stage:sd}&\frac{g^{(\ell)}_h-g^n_h}{\Delta t}+ \frac{c/\ssig}{\reps} \,\sum_{j=0}^{\ell-1} \tilde{a}_{\ell,j} \,(\mI-\Pi)\mD_h(g^{(j)}_h;\,\vO)-\frac{c/\ssig}{\reps^2} \,\sum_{j=0}^{\ell} a_{\ell,j} \,\mD^\rho_h(\vO\,\rho^{(j)}_h)
	=\\  &\hspace{5cm}-\frac{c}{\reps^2}\,\sum_{j=0}^{\ell} a_{\ell,j} \,\rsig\,g^{(j)}_h. 	\nonumber
	\end{align}
\end{subequations}
We use the mathematical induction, first to prove the AA property for the internal stages $\ell=1,\cdots,s$, by assuming a formal $\reps$-expansion \eqref{epsexpansion} for all inner stage values of $(\rho^{(\ell)}_h, g^{(\ell)}_h, \bq^{(\ell)}_h, T^{(\ell)}_h, \Phi^{(\ell)}_h)$. Similarly, first we combine \eqref{eq:FDG:GSA:stage:sa} - \eqref{eq:FDG:GSA:stage:sc} to get
\begin{align}
\label{eq:s1}
&\frac{\rho^{(\ell)}_h-\rho^n_h}{\Delta t} +\frac{c}{\ssig}\,\sum_{j=0}^{\ell-1} \tilde{a}_{\ell,j} \,\mD^g_h(\langle \vO\,g^{(j)}_h \rangle) + c\,\frac{\omega}{3\sigma_0} \,\sum_{j=0}^{\ell-1} \tilde{a}_{\ell,j} \,\mD^{div}_h\mD^{grad}(\rho^{(j)}_h) = \\ &\hspace{5cm} -\tC \,\frac{T^{(\ell)}_h-T^n_h}{\Delta t}+c\,\frac{\omega}{3\sigma_0} \,\sum_{j=0}^{\ell} \,a_{\ell,j}\,\mD^{div}_h\mD^{grad}(\rho^{(j)}_h). \nonumber
\end{align}
For the first stage $\ell = 0$, we have $\rho^{(0)}_{h}=\rho^n_h,\, g^{(0)}_{h}=g^n_h,\,\bq^{(0)}_{h}=\bq^n_h,\, T^{(0)}_{h}=T^n_h$ and $\Phi^{(0)}_{h}=\Phi^n_h$. Then for $\ell=1$, it is the same AP analysis 
for the DG-IMEX1 scheme \eqref{eq:DG}, we have
\[
\rho^{(1)}_{0,h}=\Phi^{(1)}_{0,h}+\mO(\reps)=ac(T^{(1)}_{0,h})^4/|\vO|+\mO(\reps), \quad g^{(1)}=-\frac{\ssig}{\sigma}\mD^\rho_h(\vO\rho^{(1)}_{0,h})+\mO(\reps), \quad \Phi^{(1)}_{0,h}=ac(T^{(1)}_{0,h})^4,
\]
and $T^{(1)}_{0,h}$ satisfies
\begin{align}
\label{eq:DG:diff:s1}
&\frac{a(T^{(1)}_{0,h})^4-a(T^n_{0,h})^4}{\Delta t}+ C_v \,\frac{T^{(1)}_{0,h}-T^n_{0,h}}{\Delta t} =\tilde{a}_{10}\,\mD^{div}_h\left(\frac{1}{3\sigma}\mD^{grad}_h\left(ac(T^n_{0,h})^4\right)\right) \\ &\hspace{5cm} +\frac{\omega}{3\sigma_0} \,\mD^{div}_h\mD^{grad}_h\left(a_{11}ac(T^{(1)}_{0,h})^4-\tilde{a}_{10}ac(T^{n}_{0,h})^4\right)+\mO(\reps). \nonumber
\end{align}
Now we assume for inner stages $j=1$ up to $j=\ell-1$, there holds
\[
\rho^{(j)}_{0,h}=\Phi^{(j)}_{0,h}+\mO(\reps)=ac(T^{(j)}_{0,h})^4/|\vO|+\mO(\reps), \quad g^{(j)}=-\frac{\ssig}{\sigma}\mD^\rho_h(\vO\rho^{(j)}_{0,h})+\mO(\reps), \quad \Phi^{(j)}_{0,h}=ac(T^{(j)}_{0,h})^4,
\]
and $T^{(j)}_{0,h}$ satisfies
\begin{align}
\label{eq:DG:diff:sj}
&\frac{a(T^{(j)}_{0,h})^4-a(T^n_{0,h})^4}{\Delta t}+ C_v \,\frac{T^{(j)}_{0,h}-T^n_{0,h}}{\Delta t} =\sum_{j=0}^{\ell-2} \tilde{a}_{\ell-1,j}\,\mD^{div}_h\left(\frac{1}{3\sigma}\mD^{grad}_h\left(ac(T^{(j)}_{0,h})^4\right)\right) \\ &\hspace{4cm} +\frac{\omega}{3\sigma_0} \,\mD^{div}_h\mD^{grad}_h\left(\sum_{j=0}^{\ell-1} a_{\ell-1,j}ac(T^{(j)}_{0,h})^4-\sum_{j=0}^{\ell-2} \tilde{a}_{\ell-1,j}ac(T^{(j)}_{0,h})^4\right)+\mO(\reps). \nonumber
\end{align}
By induction, we prove them also holding for the next stage $\ell$. First from \eqref{eq:FDG:GSA:stage:sb} and \eqref{eq:FDG:GSA:stage:sd}, by the definition of $\Phi=acT^4/|\vO|$, we easily get
\beq
\label{eq:sl}
\rho^{(\ell)}_{0,h}=\Phi^{(\ell)}_{0,h}+\mO(\reps)=ac(T^{(\ell)}_{0,h})^4/|\vO|+\mO(\reps), \quad g^{(\ell)}=-\frac{\ssig}{\sigma}\mD^\rho_h(\vO\rho^{(\ell)}_{0,h})+\mO(\reps), \quad \Phi^{(\ell)}_{0,h}=ac(T^{(\ell)}_{0,h})^4.
\eeq
Substituting them with \eqref{eq:FDG:GSA:stage:sc} into \eqref{eq:s1}, we get
\begin{align}
\label{eq:DG:diff:sl}
&\frac{a(T^{(\ell)}_{0,h})^4-a(T^n_{0,h})^4}{\Delta t}+ C_v \,\frac{T^{(\ell)}_{0,h}-T^n_{0,h}}{\Delta t} =\sum_{j=0}^{\ell-1} \tilde{a}_{\ell,j}\,\mD^{div}_h\left(\frac{1}{3\sigma}\mD^{grad}_h\left(ac(T^{(j)}_{0,h})^4\right)\right) \\ &\hspace{4cm} +\frac{\omega}{3\sigma_0} \,\mD^{div}_h\mD^{grad}_h\left(\sum_{j=0}^{\ell} a_{\ell,j}ac(T^{(j)}_{0,h})^4-\sum_{j=0}^{\ell-1} \tilde{a}_{\ell,j}ac(T^{(j)}_{0,h})^4\right)+\mO(\reps). \nonumber
\end{align}
By letting $\ell=s$, we have proved that \eqref{eq:sl} and \eqref{eq:DG:diff:sl} hold for all inner stages from $\ell=0$ to $s$.

For updating the solution at the time level $t^{n+1}$, \eqref{eq:FDG:GSA:n+1} is used, which in the strong form can be written as
\begin{subequations}
	\label{eq:FDG:GSA:n+1:s}
	\begin{align}
	\label{eq:FDG:GSA:n+1:sa}&\frac{\rho^{n+1}_h-\rho^n_h}{\Delta t}+\tC \,\frac{T^{n+1}_h-T^n_h}{\Delta t} = \\ & \hspace{1cm} +\frac{c}{\ssig}\,\sum_{\ell=0}^s \tilde{b}_\ell \,\mD^g_h(\langle \vO\,g^{(\ell)}_h \rangle) + c\,\frac{\omega}{3\sigma_0} \,\sum_{\ell=0}^s \tilde{b}_\ell \,\mD^{div}_h(\bq^{(\ell)}_h,\testR)+c\,\frac{\omega}{3\sigma_0} \,\sum_{\ell=0}^s \,b_\ell\,\mD^{div}_h(\bq^{(\ell)}_h,\testR), \nonumber \\ \, \nonumber \\
	\label{eq:FDG:GSA:n+1:sb}&\reps^2\,\tC \,\frac{T^{n+1}_h-T^n_h}{\Delta t} - b_s\,\rsig \,(\rho^{n+1}_h-\Phi^{n+1}_h) = \sum_{\ell=0}^{s-1} \,b_\ell\,\rsig \,(\rho^{(\ell)}_h-\Phi^{(\ell)}_h). 
	\end{align}
\end{subequations}
Submitting \eqref{epsexpansion} and \eqref{eq:sl} into \eqref{eq:FDG:GSA:n+1:s}, first from \eqref{eq:FDG:GSA:n+1:sb}, we get 
\beq
\rho^{n+1}_{0,h} = \Phi^{n+1}_{0,h}+\mO(\reps),
\eeq
then with \eqref{eq:FDG:GSA:n+1:sa}, we obtain
\begin{align}
\label{eq:DG:diff:n+1}
&\frac{a(T^{n+1}_{0,h})^4-a(T^n_{0,h})^4}{\Delta t}+ C_v \,\frac{T^{n+1}_{0,h}-T^n_{0,h}}{\Delta t} =\sum_{\ell=0}^{s} \tilde{b}_{\ell}\,\mD^{div}_h\left(\frac{1}{3\sigma}\mD^{grad}_h\left(ac(T^{(\ell)}_{0,h})^4\right)\right) \\ &\hspace{3.4cm} +\frac{\omega}{3\sigma_0} \,\mD^{div}_h\mD^{grad}_h\left(\sum_{\ell=0}^{s} b_{\ell}ac(T^{(\ell)}_{0,h})^4-\sum_{\ell=0}^{s} \tilde{b}_{\ell}ac(T^{(\ell)}_{0,h})^4\right)+\mO(\reps). \nonumber
\end{align}

\eqref{eq:DG:diff:n+1} with inner stages \eqref{eq:DG:diff:sl}, as $\reps\ra0$, becomes a high order IMEX RK discretization for the limiting diffusive equation \eqref{diffeq} and the temporal order $p$ from the high order IMEX RK scheme is maintained, namely, our DG-IMEXp scheme is AA. 

%
\section{Numerical examples}
\label{sec4}
\setcounter{equation}{0}
\setcounter{figure}{0}
\setcounter{table}{0}

In this section, some numerical tests will be performed to validate the high order accuracy, AP and AA properties of our proposed scheme. In the angular direction, the discrete coordinate $S_N$ method is used. In 1D for $\mu\in[-1,1]$, the $8$-point
Gauss quadrature rule is used. while in 2D, we take the Gauss-Chebyshev quadrature rule, with $8$ Gauss quadrature points and $4$ Chebyshev quadrature points. In space, we use the nodal DG scheme \cite{hesthaven2008nodal} with the Larangian bases, and $k$ Gaussian quadrature points are used in each space direction for a $k$-th order method. In time, we take the 3rd order type ARS(4,4,3) IMEX RK scheme \cite{ascher1997implicit}, with the double Butcher table given by
\begin{displaymath}
\begin{array}{c|c c c c c}
0 & 0&0&0& 0 & 0\\
1/2 &1/2&0&0& 0 &0\\
2/3 &11/18&1/18&0& 0 &0\\
1/2 &5/6&-5/6&1/2& 0 &0\\
1 &1/4&7/4&3/4& -7/4 &0\\
\hline
&1/4&7/4&3/4& -7/4 &0\\
\end{array}, \ \ \ \ \
\begin{array}{c|c c c c c}
0 & 0&0&0& 0 & 0\\
1/2 &0&1/2&0&0 &0\\
2/3 &0&1/6&1/2& 0 &0\\
1/2 &0&-1/2&1/2& 1/2 &0\\
1 &0&3/2&-3/2& 1/2 &1/2\\
\hline
&0&3/2&-3/2& 1/2 &1/2\\
\end{array}.
\end{displaymath}
The time step is chosen to be $\Delta t = 0.01\,h$ for the 3rd order scheme. In the following, $C_v$ and $\eps$ instead of $\tC$ and $\reps$ are specified.

In the computation, the unit of the length is taken to be centimeter ($cm$), the mass unit is gramme ($g$), the time unit is nanosecond ($ns$), the temperature unit is kilo electron-volt ($keV$) and the energy unit is $10^9$ Joules ($GJ$). Under these units, the light speed $c$ is $29.98\,cm/ns$ and the radiation constant $a$ is $0.01372\,GJ/cm^3/keV^4$, unless otherwise specified.

\subsection{Accuracy test in 1D}
\label{ex1}
We first consider a 1D example with smooth initial conditions at the equilibrium, which are given by
\beq
\label{ini1}
\rho(0,x)=(b_0+b_1\,\sin(x))^4, \quad T(0,x)=b_0+b_1\,\sin(x), \quad g(0,x,\mu)=-4\mu b_1\,\cos(x)(b_0+b_1\,\sin(x))^6,  
\eeq
where $g=-\mu\rho_x/\sigma$ and $\sigma=1/T^3$. For simplicity, all constant parameters, such as $a$, $c$ and $C_v$, are all taken to be $1$. We choose $b_0=0.8$ and $b_1=0.1$ in this test. $\sigma_0$ is taken to be $1$. Periodic boundary condition is used.

We run the solution with 3rd order method in time. While in space, the nodal DG scheme with $k=2$ to $4$ Gaussian quadrature points are considerred, resulting 2nd, 3rd or 4th order respectively. $\eps=1, 10^{-2}, 10^{-6}$ in three different regimes are taken. Since the exact solutions are not known, we compute the numerical errors by comparing the solutions on two consecutive mesh sizes. In Table \ref{tab1}, we show the errors and orders in $L^1$ and $L^\infty$ norms for different $k$'s and different $\eps$'s. Due to the small time step condition, the spatial errors are dominate, and we can observe $k$th order for the three $\eps$'s, which verfies that our scheme is asymptotic preserving and asymptotically accurate.

\begin{table}[htbp]
\caption{Numerical errors and orders of $L^1$ and $L^\infty$ norms for accuracy test in 1D. $T=0.2$. $k$ Gaussian quadrature points for $k=2,3,4$ in space are used.}
\begin{center}
\begin{tabular}{|c|c|c|c|c|c|c|}\hline
 \multicolumn{2}{|c|}{}  &   N  &  $L^1$ error &   order &  $L^\infty$ error & order \\ \hline
  \multirow{12}{1cm}{$k=2$}  &  \multirow{4}{1.5cm}{$\eps=1$}
   &    20 &     2.02E-03 &       --&     7.39E-03 &       -- \\ \cline{3-7}
&  &    40 &     6.36E-04 &     1.67&     2.80E-03 &     1.40 \\ \cline{3-7}
&  &    80 &     1.81E-04 &     1.82&     8.22E-04 &     1.77 \\ \cline{3-7}
&  &   160 &     4.34E-05 &     2.06&     1.91E-04 &     2.10 \\ \cline{2-7}
   &  \multirow{4}{1.5cm}{$\eps=10^{-2}$}
   &    20 &     2.29E-03 &       --&     9.63E-03 &       -- \\ \cline{3-7}
&  &    40 &     6.24E-04 &     1.88&     2.40E-03 &     2.01 \\ \cline{3-7}
&  &    80 &     1.55E-04 &     2.01&     6.03E-04 &     1.99 \\ \cline{3-7}
&  &   160 &     3.86E-05 &     2.00&     1.51E-04 &     2.00 \\ \cline{2-7}
   &  \multirow{4}{1.5cm}{$\eps=10^{-6}$}
   &    20 &     2.31E-03 &       --&     9.72E-03 &       -- \\ \cline{3-7}
&  &    40 &     6.23E-04 &     1.89&     2.39E-03 &     2.02 \\ \cline{3-7}
&  &    80 &     1.55E-04 &     2.01&     6.02E-04 &     1.99 \\ \cline{3-7}
&  &   160 &     3.86E-05 &     2.00&     1.51E-04 &     2.00 \\ \hline
 
   \multirow{12}{1cm}{$k=3$} & \multirow{4}{1.5cm}{$\eps=1$}
   &   20 &     1.28E-04 &       --&     4.96E-04 &       -- \\ \cline{3-7}
&  &   40 &     1.98E-05 &     2.69&     8.07E-05 &     2.62 \\ \cline{3-7}
&  &   80 &     2.58E-06 &     2.94&     1.14E-05 &     2.83 \\ \cline{3-7}
&  &  160 &     3.34E-07 &     2.95&     1.51E-06 &     2.91 \\ \cline{2-7}
   &  \multirow{4}{1.5cm}{$\eps=10^{-2}$}
   &   20 &     1.23E-04 &       --&     5.09E-04 &       -- \\ \cline{3-7}
&  &   40 &     1.64E-05 &     2.91&     7.20E-05 &     2.82 \\ \cline{3-7}
&  &   80 &     2.03E-06 &     3.01&     8.96E-06 &     3.01 \\ \cline{3-7}
&  &  160 &     2.53E-07 &     3.00&     1.11E-06 &     3.01 \\ \cline{2-7}
   &  \multirow{4}{1.5cm}{$\eps=10^{-6}$}
   &   20 &     1.23E-04 &       --&     5.08E-04 &       -- \\ \cline{3-7}
&  &   40 &     1.64E-05 &     2.91&     7.18E-05 &     2.82 \\ \cline{3-7}
&  &   80 &     2.03E-06 &     3.01&     8.93E-06 &     3.01 \\ \cline{3-7}
&  &  160 &     2.53E-07 &     3.00&     1.11E-06 &     3.01 \\ \cline{1-7}

   \multirow{12}{1cm}{$k=4$} & \multirow{4}{1.5cm}{$\eps=1$}
   &    20 &     9.84E-06 &       --&     4.85E-05 &       -- \\ \cline{3-7}
&  &    40 &     7.52E-07 &     3.71&     4.25E-06 &     3.51 \\ \cline{3-7}
&  &    80 &     4.37E-08 &     4.10&     2.51E-07 &     4.08 \\ \cline{3-7}
&  &   160 &     2.71E-09 &     4.01&     1.56E-08 &     4.01 \\ \cline{2-7}
&  \multirow{4}{1.5cm}{$\eps=10^{-2}$}
   &    20 &     7.18E-06 &       --&     3.35E-05 &       -- \\ \cline{3-7}
&  &    40 &     4.37E-07 &     4.04&     1.99E-06 &     4.07 \\ \cline{3-7}
&  &    80 &     2.72E-08 &     4.01&     1.27E-07 &     3.98 \\ \cline{3-7}
&  &   160 &     1.70E-09 &     4.00&     7.93E-09 &     4.00 \\ \cline{2-7}
&  \multirow{4}{1.5cm}{$\eps=10^{-6}$}
   &    20 &     7.27E-06 &       --&     3.42E-05 &       -- \\ \cline{3-7}
&  &    40 &     4.37E-07 &     4.06&     1.98E-06 &     4.11 \\ \cline{3-7}
&  &    80 &     2.72E-08 &     4.01&     1.26E-07 &     3.97 \\ \cline{3-7}
&  &   160 &     1.86E-09 &     3.87&     7.91E-09 &     4.00 \\ \hline
\end{tabular}
\end{center}
\label{tab1}
\end{table}

\subsection{Marshak wave}
\label{ex2}
The second example is the Marshak wave problem in 1D \cite{Sun2015an}. We first take the absorption/emission coefficient to be $\sigma=300/T^3/cm$ for the Marshak wave-2B problem. The specific heat is $0.1\,GJ/g/keV$, the density is $3.0\,g/cm^3$ and $\eps=1$. The initial temperature $T$ is set to be $10^{-6}\,keV$ everywhere. A constant isotropic incident radiation intensity with a Planckian distribution at $1\,keV$ is set at the left boundary. For this problem, the inflow-outflow close-loop boundary condition \cite{Jang2015high,Peng2020stability} is used. 

\begin{figure}
	\centering
	\includegraphics[width=3.2in,clip]{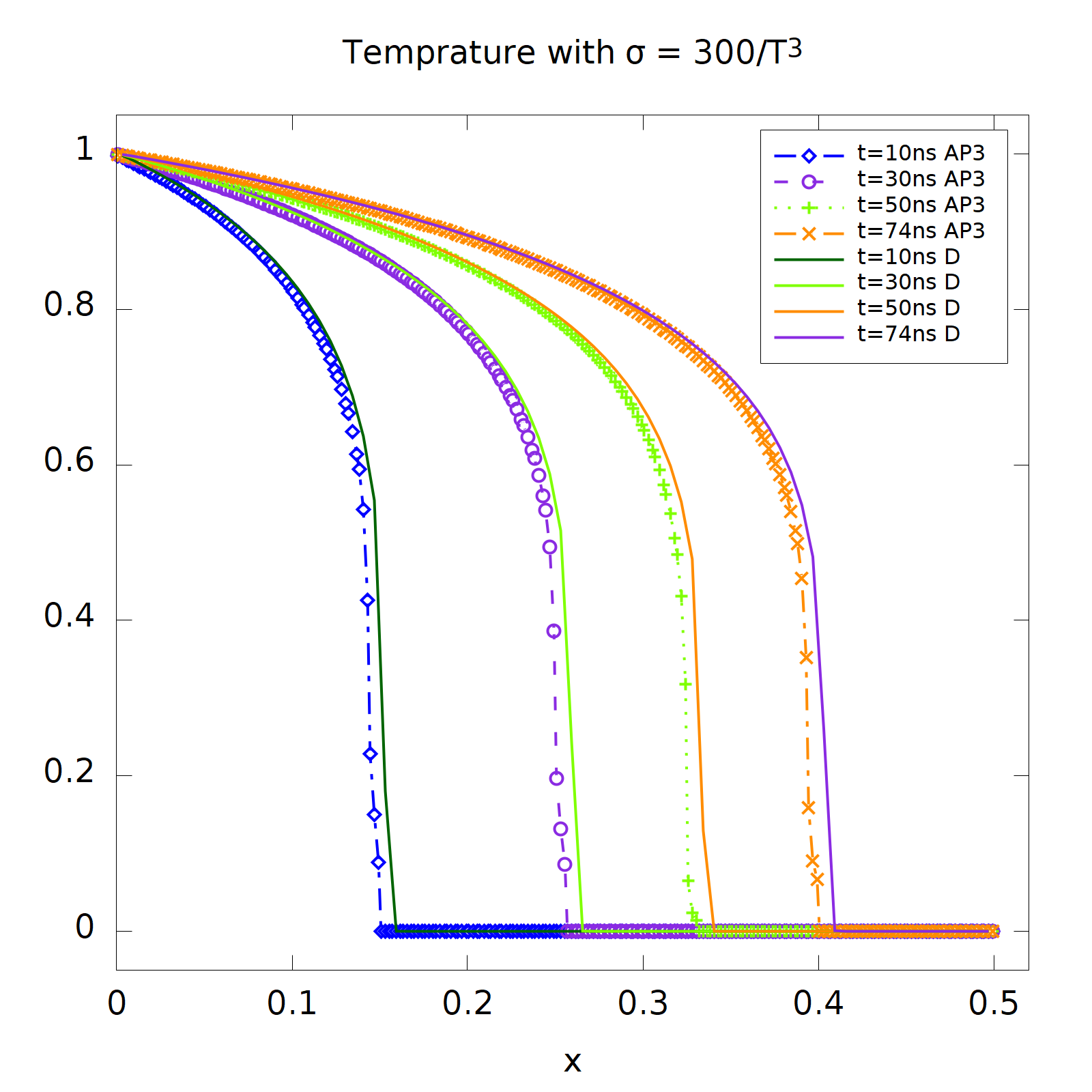}
	\caption{Marshak wave-2B problem. 3rd order scheme with $N=80$. }
	\label{fig1}
\end{figure}

We take the 3rd method both in space and in time with $N=80$. For the 3rd order nodal DG scheme, a double minmod slope limiter \cite{McClarren2008} is used to control the numerical oscillations at the propagation front. For the Marshak wave-2B problem, we take $\sigma_0=300$ as the reference opacity. Due to $\sigma_0$ is relatively large, the solution to the GRTEs would be very close to the solution for the limiting diffusive equation \eqref{diffeq}.
In Fig. \ref{fig1}, we show the solutions obtained from our 3rd order AP scheme, which is denoted as ``AP3'' and compared to the limiting scheme 
for the diffusive equation \eqref{diffeq} (directly implemented for \eqref{diffeq}) denoted as ``D''. We can observe that the two sets of solutions are very close at $t=10, 30, 50, 74$ nano seconds. The front positions also match those in \cite{Sun2015an} (Fig. 4 \& 5), which show that our scheme has the right diffusive limit.  

\begin{figure}
	\centering
	\includegraphics[width=3.2in,clip]{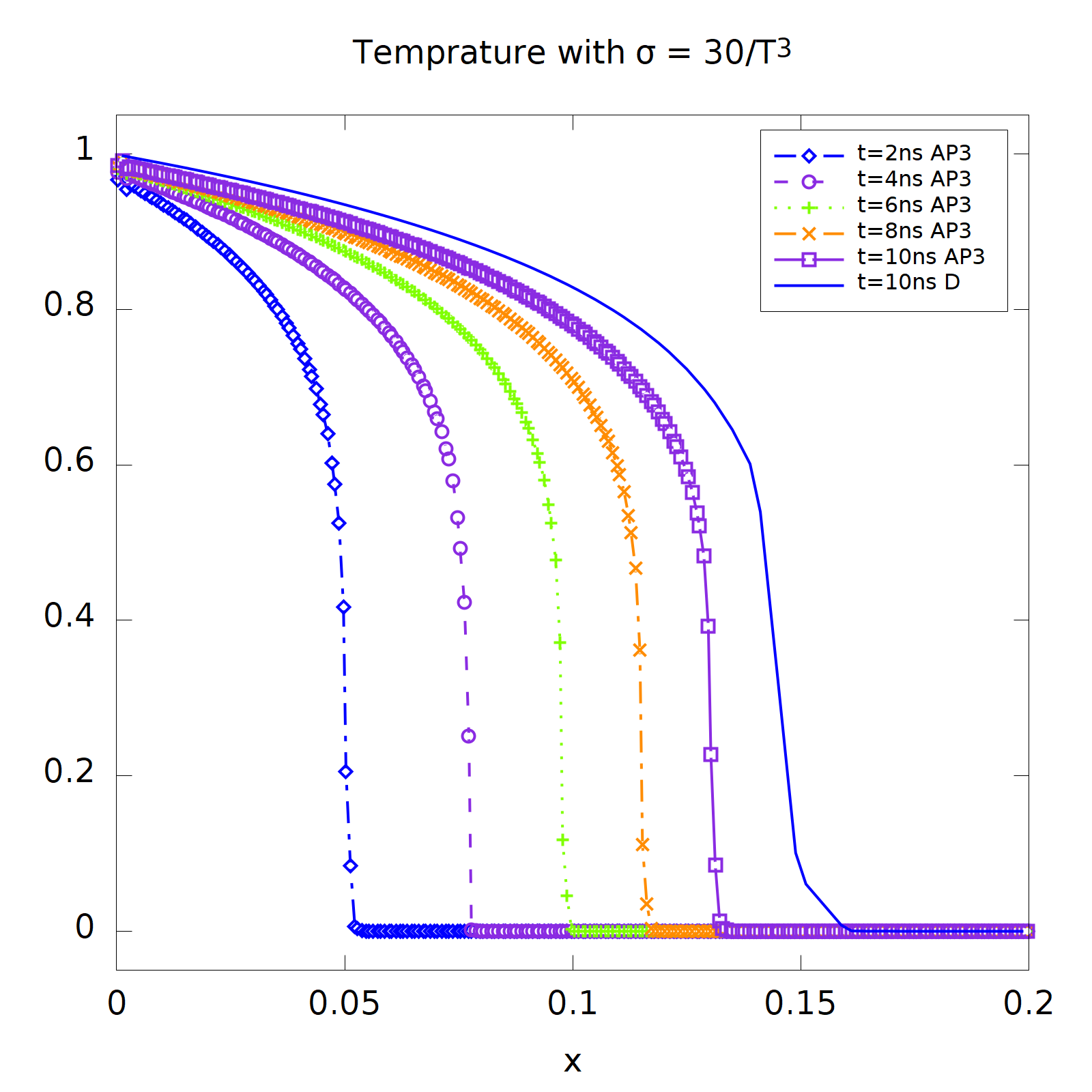}
	\caption{Marshak wave-2A problem. 3rd order scheme with $N=80$.. }
	\label{fig2}
\end{figure}

Then we take the absorption/emission coefficient is taken to be $\sigma=30/T^3/cm$ for the Marshak wave-2A problem. Other settings are the same as the Marshak wave-2B problem and the same boundary condition is used. Here $\sigma_0=30$. For this case, we should observe that the solution of the AP scheme would deviate from the solution in the diffusive limit. In Fig. \ref{fig2}, we show the solution from the AP3 scheme at $t=0.2, 0.4, 0.6, 0.8$ nano seconds, and compare them to the solution at $t=1ns$ from the D scheme. As we can see that, the solution of AP3 scheme propagates slower than that from the D scheme. Also our solutions have comparable propagation fronts as those in \cite{Sun2015an} (Fig. 6 \& 7), including the solution from the D scheme.

\subsection{Accuracy test in 2D}
\label{ex4}
Now we test the order of convergence in the 2D case. We consider smooth initial conditions at the equilibrium similar to the 1D problem, which are given by
\beq
\rho(0,x,y) = \left((a1+b1\sin(x))(a2+b2\sin(y)\right)^4, \quad
T(0,x,y)=(a1+b1\sin(x))(a2+b2\sin(y)),
\eeq
and $g(0,x,y)=-\vec{\Omega}\cdot\nabla \rho(0,x,y)/\sigma$. Here $\sigma=1$. Similarly all constant parameters, such as $a$, $c$ and $C_v$, are all taken to be $1$. We choose $a1=a2=0.8$ and $b1=b2=0.1$ in this test. Periodic boundary conditions are used in both space directions.

We run the solution with 3rd order method in time. While in space, the nodal DG scheme with $k=2$ to $4$ Gaussian quadrature points in the two space directions are considerred, resulting 2nd, 3rd or 4th order respectively. $\eps=1, 10^{-2}, 10^{-6}$ in three different regimes are taken. Similarly we compute the numerical errors by comparing the solutions on two consecutive mesh sizes. In Table \ref{tab2}, we show the errors and orders in $L^1$ and $L^\infty$ norms for different $k$'s and different $\eps$'s. Similar results as the 1D problem are obtained.

\begin{table}[htbp]
\caption{Numerical errors and orders of $L^1$ and $L^\infty$ norms for accuracy test in 2D. $T=0.01$. $k$ Gaussian quadrature points for $k=2,3,4$ in both two space directions are used.}
\begin{center}
\begin{tabular}{|c|c|c|c|c|c|c|}\hline
\multicolumn{2}{|c|}{} &  $N \times N$  &  $L^1$ error &   order &  $L^\infty$ error & order \\ \hline
     \multirow{9}{1cm}{$k=2$} & \multirow{3}{1.5cm}{$\eps=1$}
  &  $16\times 16$ &     1.63E-03 &       --&     1.55E-02 &       -- \\ \cline{3-7}
& &  $32\times 32$ &     4.15E-04 &     1.97&     4.57E-03 &     1.77 \\ \cline{3-7}
& &  $64\times 64$ &     1.05E-04 &     1.99&     1.19E-03 &     1.94 \\ \cline{2-7}
&	 \multirow{3}{1.5cm}{$\eps=10^{-2}$}
  &  $16\times 16$ &     1.65E-03 &       --&     1.63E-02 &       -- \\ \cline{3-7}
& &  $32\times 32$ &     4.61E-04 &     1.84&     5.24E-03 &     1.64 \\ \cline{3-7}
& &  $64\times 64$ &     1.49E-04 &     1.63&     1.79E-03 &     1.55 \\ \cline{2-7}
&  \multirow{3}{1.5cm}{$\eps=10^{-6}$}
  &  $16\times 16$ &     1.66E-03 &       --&     1.64E-02 &       -- \\ \cline{3-7}
& &  $32\times 32$ &     4.68E-04 &     1.82&     5.32E-03 &     1.62 \\ \cline{3-7}
& &  $64\times 64$ &     1.52E-04 &     1.62&     1.83E-03 &     1.54 \\ \hline

     \multirow{9}{1cm}{$k=3$} & \multirow{3}{1.5cm}{$\eps=1$}
  &  $16\times 16$ &     1.16E-04 &       --&     1.03E-03 &       -- \\ \cline{3-7}
& &  $32\times 32$ &     1.45E-05 &     3.00&     1.49E-04 &     2.80 \\ \cline{3-7}
& &  $64\times 64$ &     1.86E-06 &     2.96&     1.92E-05 &     2.95 \\ \cline{2-7}
  &	 \multirow{3}{1.5cm}{$\eps=10^{-2}$}
  &  $16\times 16$ &     1.27E-04 &       --&     1.12E-03 &       -- \\ \cline{3-7}
& &  $32\times 32$ &     1.94E-05 &     2.71&     1.99E-04 &     2.49 \\ \cline{3-7}
& &  $64\times 64$ &     2.74E-06 &     2.83&     2.98E-05 &     2.74 \\ \cline{2-7}
  &  \multirow{3}{1.5cm}{$\eps=10^{-6}$}
  &  $16\times 16$ &     1.29E-04 &       --&     1.13E-03 &       -- \\ \cline{3-7}
& &  $32\times 32$ &     1.97E-05 &     2.71&     2.02E-04 &     2.48 \\ \cline{3-7}
& &  $64\times 64$ &     2.80E-06 &     2.82&     2.95E-05 &     2.78 \\ \hline

     \multirow{9}{1cm}{$k=4$} & \multirow{3}{1.5cm}{$\eps=1$}
  &  $16\times 16$ &     8.88E-06 &       --&     9.11E-05 &       -- \\ \cline{3-7}
& &  $32\times 32$ &     5.77E-07 &     3.95&     7.46E-06 &     3.61 \\ \cline{3-7}
& &  $64\times 64$ &     3.94E-08 &     3.87&     5.11E-07 &     3.87 \\ \cline{2-7}
&	 \multirow{3}{1.5cm}{$\eps=10^{-2}$}
  &  $16\times 16$ &     1.04E-05 &       --&     1.17E-04 &       -- \\ \cline{3-7}
& &  $32\times 32$ &     6.95E-07 &     3.90&     1.01E-05 &     3.53 \\ \cline{3-7}
& &  $64\times 64$ &     5.43E-08 &     3.68&     7.34E-07 &     3.79 \\ \cline{2-7}
     &  \multirow{3}{1.5cm}{$\eps=10^{-6}$}
  &  $16\times 16$ &     1.05E-05 &       --&     1.19E-04 &       -- \\ \cline{3-7}
& &  $32\times 32$ &     6.99E-07 &     3.91&     1.02E-05 &     3.54 \\ \cline{3-7}
& &  $64\times 64$ &     4.72E-08 &     3.89&     7.32E-07 &     3.80 \\ \hline

		\end{tabular}
	\end{center}
	\label{tab2}
\end{table}
 
\subsection{Tophat test}
\label{ex5}
This is a 2D gray radiative transfer test problem, which has already been studied in \cite{Gentile2001,Sun2015an, Hammer2019, Shi2018}. The computational domain is $[0,7]\times[-2,2]$. The dense and opaque material with density $10\,g/cm^3$ and opacity $\sigma=2000\,cm^{-1}$ is located in the following regions: $(3,4)\times(-1,1), (0,2.5)\times(-2,-0.5), (0,2.5)\times(0.5,2),(4.5,7)\times(-2,-0.5), (4.5,7)\times(0.5,2),(2.5,4.5)\times(-2,-1.5),(2.5,4.5)\times(1.5,2)$. The pipe, which has density $0.01\,g/cm^3$ and opacity $\sigma=0.2\,cm^{-1}$, occupies all other regions. The heat capacity is $0.1\,GJ/g/keV$ and $\eps=1$. 

\begin{figure}
	\centering
	\includegraphics[width=3.2in,clip]{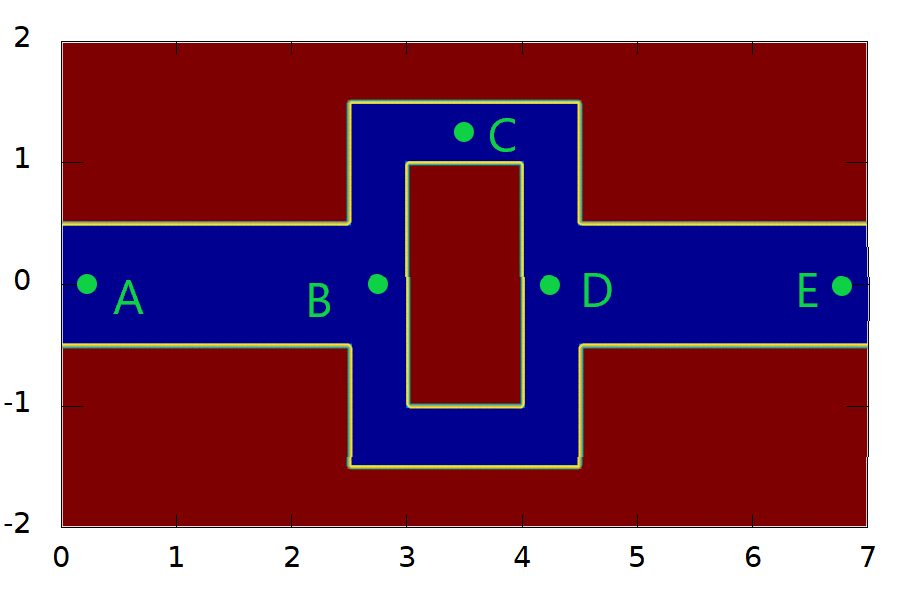}
	\caption{Initial configuration for the Tophat test. In the pipe, $\rho=0.01\,g/cm^3$ and  $\sigma=0.2\,cm^{-1}$, and outer $\rho=10\,g/cm^3$ and opacity $\sigma=2000\,cm^{-1}$. Five probes are placed to monitor the change of temperature.}
	\label{fig3}
\end{figure}

Initially, the material has a temperature $0.05\,keV$ everywhere, and the radiation and material temperature are in equilibrium. A heating source with a fixed temperature $0.5\,keV$ is located on the left boundary for $-0.5<y<0.5$. Others are outflow boundary conditions.
Five probes (A,B,C,D,E) are placed at $(0.25,0), (2.75,0), (3.5,1.25), (4.25,0), (6.75,0)$ to monitor the change of temperature in the thin opacity material. See Fig. \ref{fig3} for the initial configuration.

For this problem, we take $\sigma_0=10000$. 3rd order in space is used. For this long time simulation, we choose the first order IMEX scheme \eqref{eq:DG} in time.
Besides, due to the rapid change of the opacity at the interface between the two materials, the simple numerical fluxes \eqref{eq:vg:upwind:L-1} and \eqref{eq:flux} cannot work properly. Instead, the opacity should be taken into account for computing the numerical fluxes at the interface. For our scheme \eqref{eq:SDG} under the micro-macro decomposition, it mainly affects the fluxes related to $\rho$. Taking \eqref{eq:flux:1} as an example, instead of purely downwind for $\hat\rho=\rho^+$, at the interface, we modify it to be $\hat \rho = (1-\omega_1) \rho^+ + \omega_1 \rho^-$, where $\omega_1=\exp(-c\rsig^+\,dt/\reps^2)$ and $\rsig^+$ is taking from the same side as $\rho^+$. This weight is the exponential function appeared in the time dependent evolution solution for the GRTEs \eqref{rte}, which is used to desigh a unified gas-kinetic scheme, e.g., see (3.6) in \cite{Sun2015an}.  

\begin{figure}
	\centering
	\includegraphics[width=3.2in,clip]{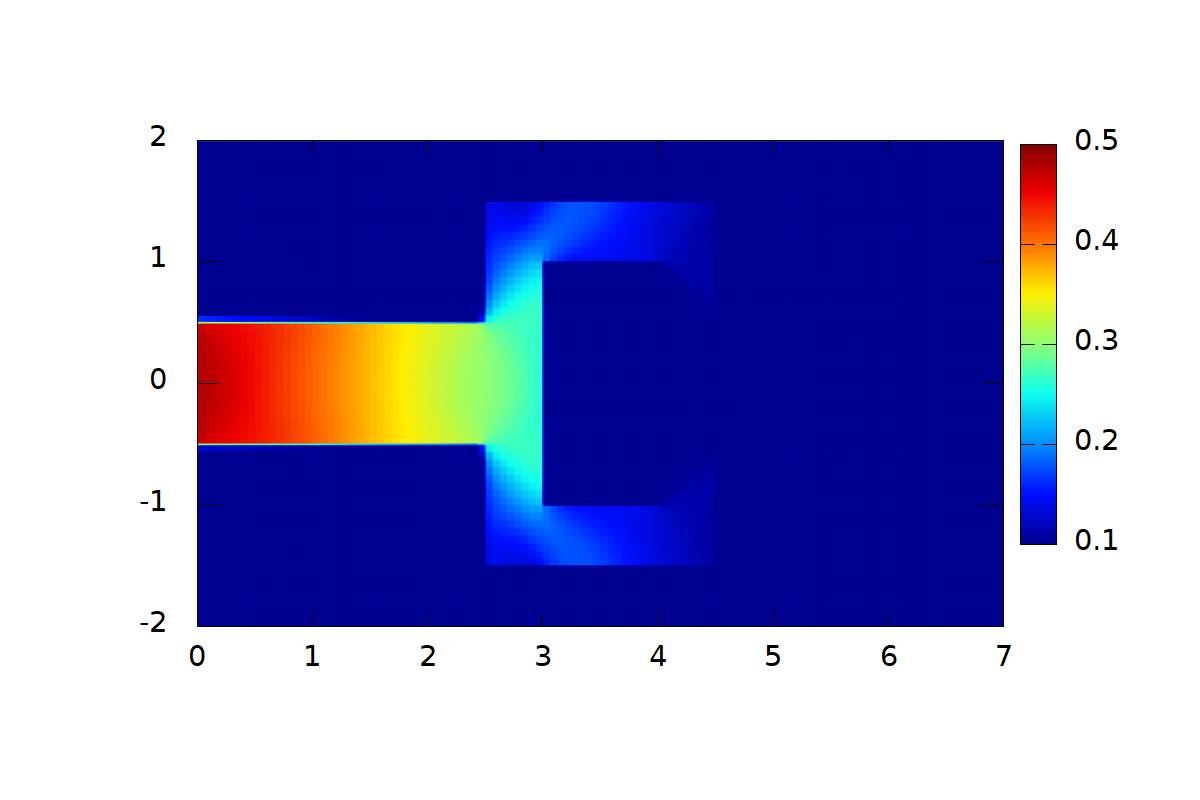}
	\includegraphics[width=3.2in,clip]{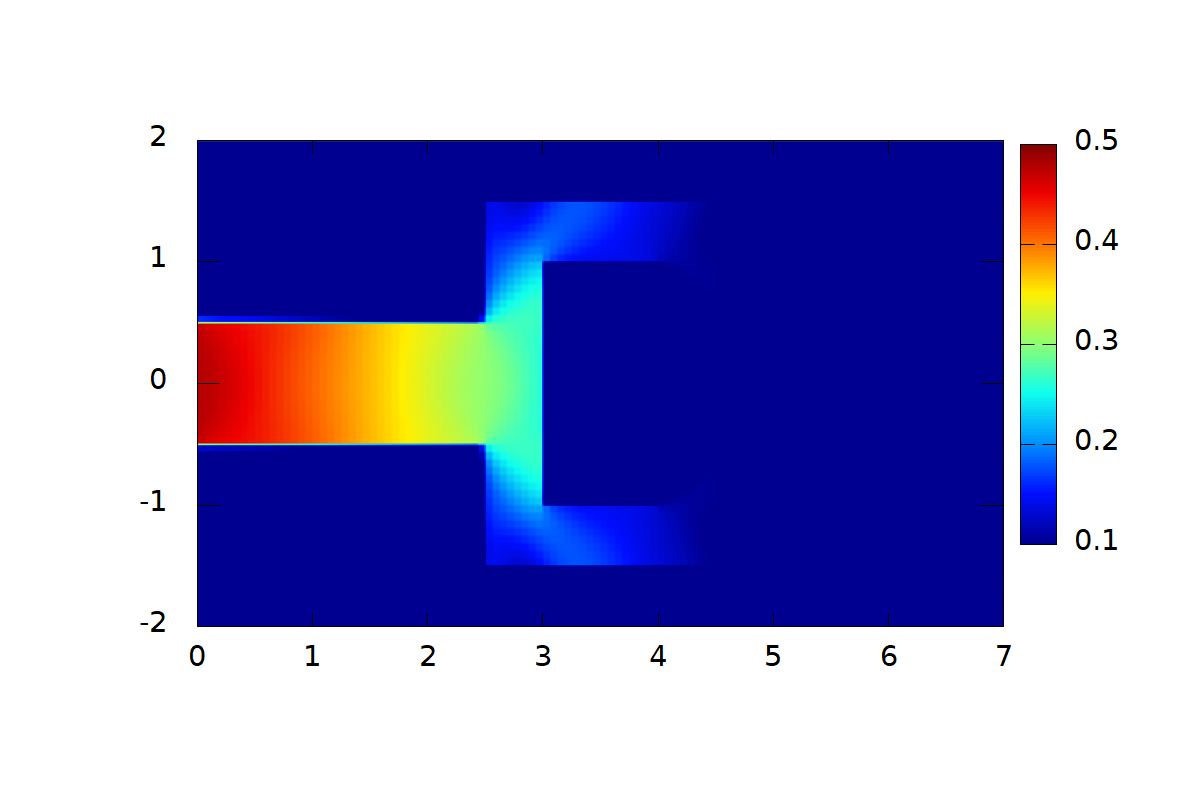}
	\caption{Tophat test. The radiative temperature $T_{r}$ and the material temperature $T$ at $10\,ns$. Mesh is $112\times64$. }
	\label{fig4}
\end{figure}

\begin{figure}
	\centering
	\includegraphics[width=3.2in,clip]{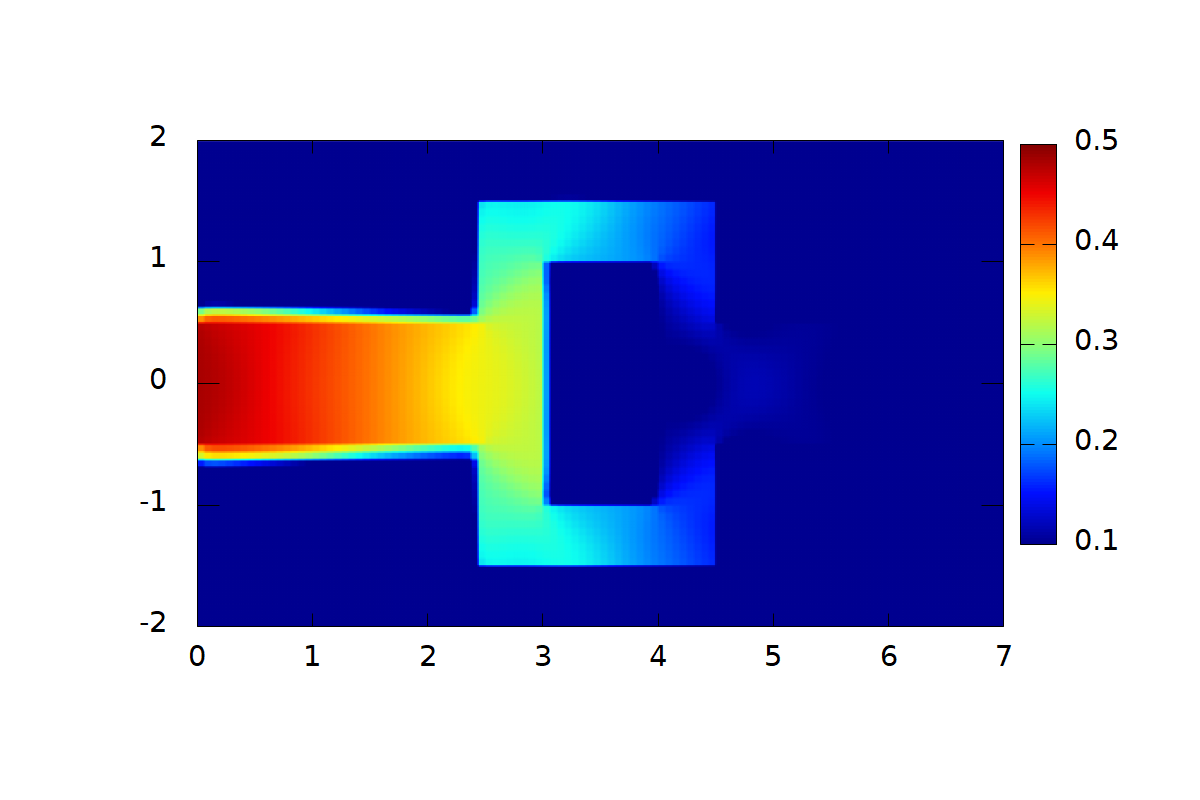}
    \includegraphics[width=3.2in,clip]{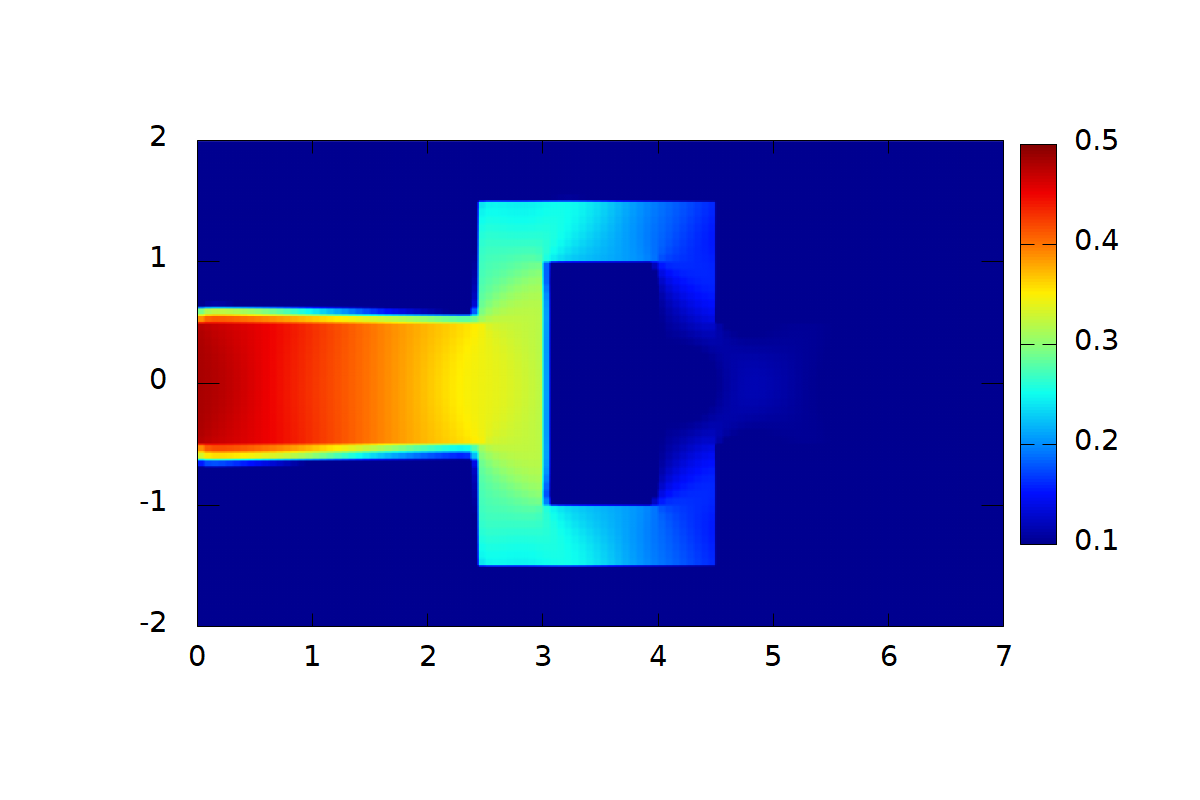}
	\caption{Tophat test. The radiative temperature $T_{r}$ and the material temperature $T$ at $100\,ns$. Mesh $112\times64$.  }
	\label{fig5}
\end{figure}

\begin{figure}
	\centering
	\includegraphics[width=2.6in,clip]{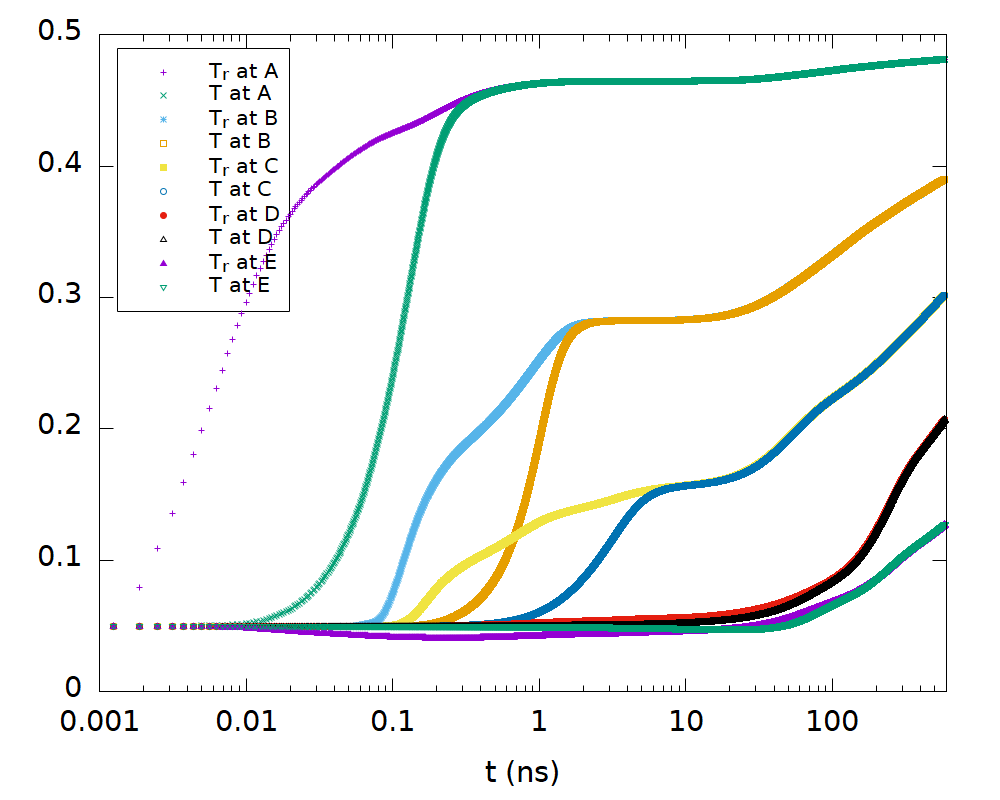}
	\caption{Tophat test. The change of the radiative temperature $T_r$ and the material temperature $T$ at the five probes upto $1000\,ns$. Mesh is $56\times32$. }
	\label{fig6}
\end{figure}

\begin{figure}
	\centering
	\includegraphics[width=2.6in,clip]{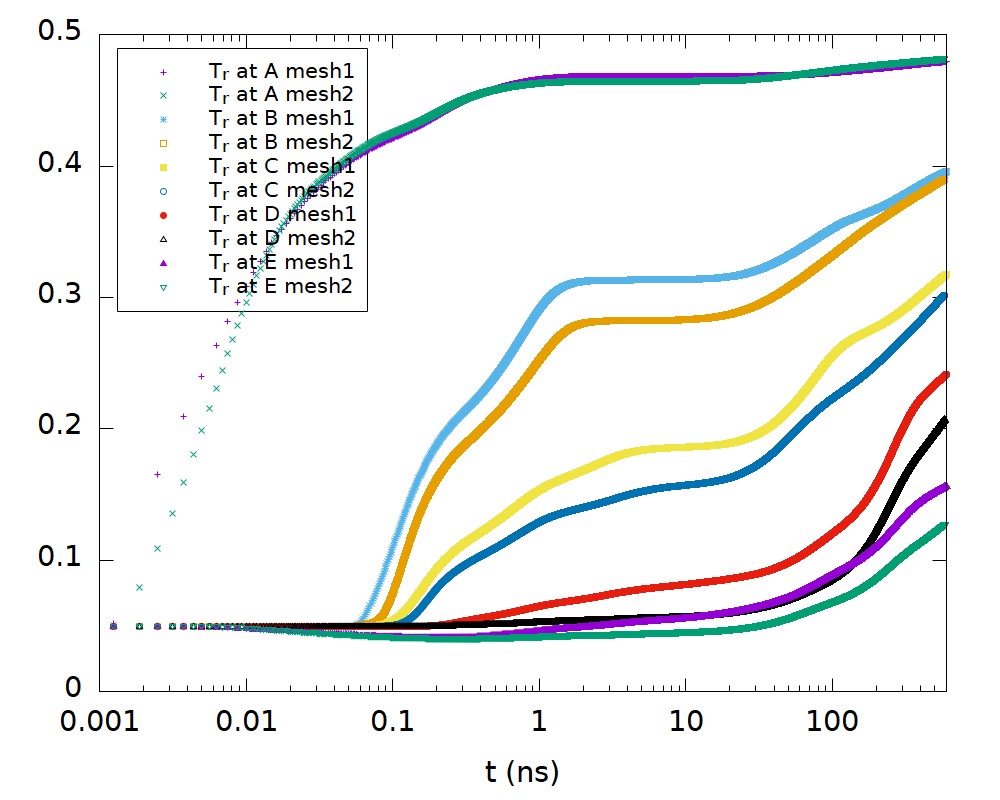}
	\includegraphics[width=2.6in,clip]{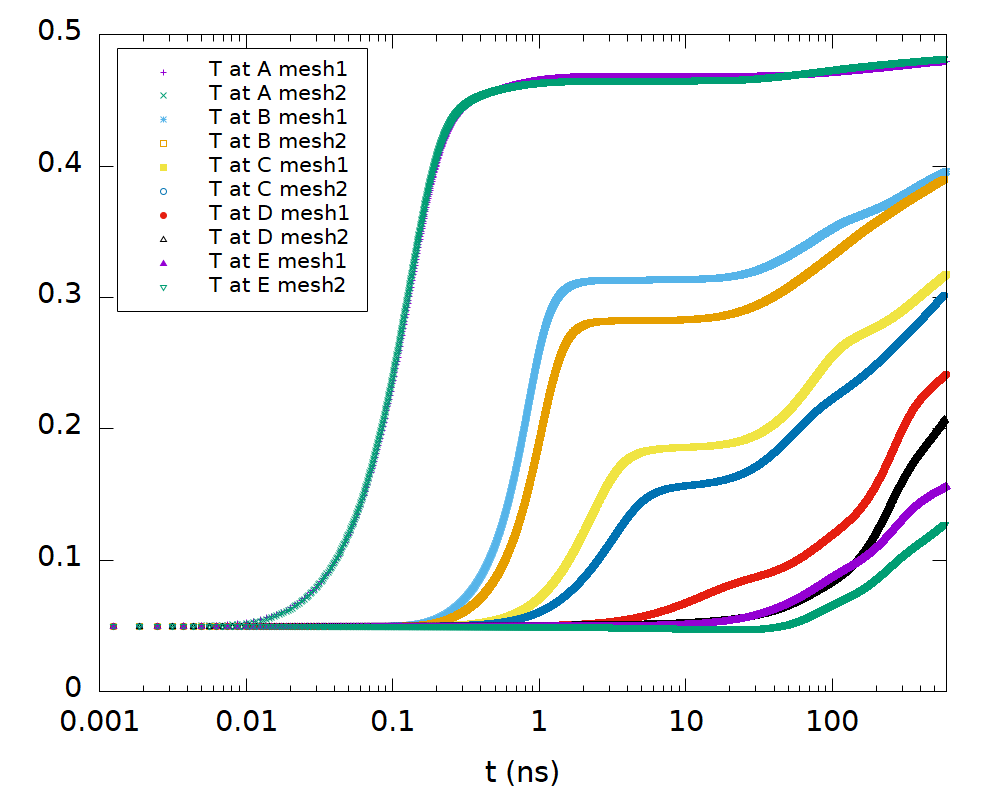}
	\caption{Tophat test. A comparison of the radiative temperature $T_r$ (left) and the material temperature $T$ (right) at the five probes upto $600\,ns$. Mesh1 is $56\times32$ and mesh2 is $112\times64$. }
	\label{fig7}
\end{figure}

We first take mesh size to be $112\times64$ and compute the solution upto $t=400ns$. In Fig. \ref{fig4} and Fig. \ref{fig5}, we show the radiative temperature, which is defined as $T_r=(|\vO|\langle I\rangle/ac)^{1/4}$, and the material temperature $T$ at $t=10$ and $100$ nano seconds, respectively. We also show the change of the radiative temperature $T_r$ and the material temperature $T$ at the five probles on the mesh $56\times42$ in Fig. \ref{fig6}. It is similar to the results in \cite{Sun2015an} (Fig. 10). We note that it is not easy to do a good comparison for the change of $T_r$ and $T$, especially when the time is not long enough, e.g., see \cite{Hammer2019} (Fig. 5). Here we compare $T_r$ and $T$ from our method on two different meshes $112\times64$ and $56\times32$, see Fig. \ref{fig7}. Relative convergent results can be seen. 

%
\section{Conclusion}
\label{sec5}
\setcounter{equation}{0}
\setcounter{figure}{0}
\setcounter{table}{0}
In this work, a class of high order asymptotic preserving DG IMEX schemes is developed for the gray radiative transfer equations. A weighted linear diffusive term is added for pernalization, so that the resulting scheme in the diffusive limit, requires only a hyperbolic type time step restriction. The scheme has been formally proved to be AP and AA for well-prepared initial conditions. Numerical results have validated the high order, AP and AA properties of the scheme, and good performance for the 1D Marshak wave and 2D Tophat test problems. The efficiency of the scheme as compared to some other schemes will be investigated in our future work.

\section*{Acknowledgement}

All the authors acknowledge support by the Science Challenge Project No. TZ2016002. T. Xiong also acknowledges support by NSFC grant No. 11971025, NSF grant of Fujian Province No. 2019J06002.

%

\bibliographystyle{siam}
\bibliography{refer}
\end{document}